\documentclass[journal,twoside]{IEEEtran}
\usepackage{amsmath,amssymb,amsfonts}
\usepackage{algorithmic}
\usepackage{algorithm}
\usepackage{array}
\usepackage{textcomp}
\usepackage{stfloats}
\usepackage{url}
\usepackage{verbatim}
\usepackage{graphicx}
\usepackage{cite}
\usepackage{flushend}
\usepackage{multicol}
\usepackage{multirow}
\usepackage{booktabs}
\usepackage{makecell}

\usepackage{xcolor}

\usepackage{enumitem}
\setenumerate[1]{itemsep=0pt,partopsep=0pt,parsep=\parskip,topsep=0pt}
\setitemize[1]{itemsep=0pt,partopsep=0pt,parsep=\parskip,topsep=0pt}
\setdescription{itemsep=0pt,partopsep=0pt,parsep=\parskip,topsep=0pt}

\usepackage{caption}
\usepackage{subcaption}

\usepackage{dsfont}

\allowdisplaybreaks[4]  

\newtheorem{assumption}{Assumption}
\newtheorem{theorem}{Theorem}
\newtheorem{lemma}{Lemma}
\newtheorem{definition}{Definition}

\begin{document}

\title{A Data-Driven Real-Time Optimal Power Flow Algorithm Using Local Feedback}

\author{Heng Liang, \IEEEmembership{Student Member, IEEE}, Yujin Huang, and Changhong Zhao, \IEEEmembership{Senior Member, IEEE}
\thanks{This work was supported by Hong Kong Research Grants Council under General Research Fund No. 14212822. Corresponding author: Changhong Zhao.}
\thanks{The authors are with the Department of Information Engineering, the Chinese University of Hong Kong, Hong Kong, China. Emails: \{lh021, hy123, chzhao\}@ie.cuhk.edu.hk}
}


\maketitle

\begin{abstract}
The increasing penetration of distributed energy resources (DERs) adds variability as well as fast control capabilities to power networks. Dispatching the DERs based on local information to provide real-time optimal network operation is the desideratum. In this paper, we propose a data-driven real-time algorithm that uses only the local measurements to solve time-varying AC optimal power flow (OPF). Specifically, we design a learnable function that takes the local feedback as input in the algorithm. The learnable function, under certain conditions, will result in a unique stationary point of the algorithm, which in turn transfers the OPF problems to be optimized over the parameters of the function. We then develop a stochastic primal-dual update to solve the variant of the OPF problems based on a deep neural network (DNN) parametrization of the learnable function, which is referred to as the training stage. We also design a gradient-free alternative to bypass the cumbersome gradient calculation of the nonlinear power flow model. The OPF solution-tracking error bound is established in the sense of universal approximation of DNN. Numerical results on the IEEE $37$-bus test feeder show that the proposed method can track the time-varying OPF solutions with higher accuracy and faster computation compared to benchmark methods.
\end{abstract}


\section{Introduction}
\label{sec:introduction}
\IEEEPARstart{O}{ptimal} power flow (OPF) is a fundamental optimization problem that determines the best operating point for dispatchable devices to achieve an optimal system-wide objective subject to the physical laws and safety constraints. Extensive algorithms have been developed to solve the important OPF problems efficiently. Most of those algorithms require intense communications to reach a converged solution and then apply the solution to the power network, and hence they are only applicable to a slow timescale operation and do not have stability guarantees from the control perspective. However, the deployment of massive distributed energy resources (DERs, such as wind and solar generations, smart appliances, energy storage devices) introduces increased variations as well as fast control capabilities to the power grids, which poses an urgent need for real-time OPF algorithms that conduct optimization and control simultaneously. Nevertheless, the DERs are widespread and distributed in nature, and may not have adequate real-time communication capabilities to participate in the system-level OPF process. It is of great interest to leverage the available local feedback to continuously drive the power network towards the time-varying OPF solutions with strong guarantees for accuracy and responsiveness. 

Numerous efforts have been focused on developing OPF algorithms that run on a faster timescale. The conventional methods, either centralized \cite{low2014convex,lavaei2012zeroduality} or distributed \cite{dall2013distributed,erseghe2014distributed,peng2016distributed}, often consider the steady-state optimization. They need to conduct numbers of iterations until the algorithms converge and then apply the solution to the power network, which are not adaptive to the future power grids as the rapid variations of DERs require solving a large number of OPF instances within a limited time frame. To address the overwhelming computations, machine learning techniques have been utilized to accelerate the OPF solvers \cite{pineda2020datadriven,biagioni2022learning}, or predict the OPF solutions \cite{pan2021dcopf,singh2022learning}. Most of them require a large training dataset with pre-solved OPF solutions to obtain a well-trained model, while acquiring such a dataset would be non-trivial and expensive. To circumvent the efforts of preparing the dataset with ground truth, reference \cite{huang2024unsupervised} proposed an unsupervised learning approach by collectively minimizing the OPF objectives and constraint violations, while \cite{gupta2022controlling} performed the training in the dual domain to adaptively handle the constraints in the OPF problems. Although the methods in \cite{huang2024unsupervised,gupta2022controlling} demonstrated appealing results, they need to collect the required information, such as the given load demands, from all buses of the power network. However, not all buses are monitored and capable of communicating in real time, which prohibits the applications of these methods in practical power grids.


The rapid fluctuation of renewable generation in power grids necessitates real-time algorithms that can perform time-varying optimization. Reference \cite{tang2017realtime} proposed a real-time algorithm based on the quasi-Newton method that incorporates the second-order information to agilely steer the output of DERs towards the fast-changing OPF solutions. Reference \cite{anese2018optimal} developed a primal-dual controller that iteratively actuates the power grids with intermediate decisions and updates the decisions based on the system feedback. References \cite{tang2018feedback,tang2022running} further established the asymptotic bound of the primal-dual controller to the (local) optimum of time-varying nonconvex OPF problems. Those methods assume the power grids are equipped with reliable real-time communication networks, which is not often satisfied in practice, especially in distribution networks. With the success of learning-based control \cite{feng2024stability,cui2023reinforcement}, some recent works utilize machine learning models to incorporate the OPF solution process into local controllers to bypass the real-time communications. Reference \cite{yuan2024learning} designed an incremental local Volt/Var controller with a learnable neural network, in which the neural network was trained in a supervised learning manner by using many instances of pre-solved OPF solutions. Reference \cite{yuan2024unsupervised} proposed an unsupervised learning approach to train a local controller, which restricts the objective as voltage deviations and ignores inequality constraints (e.g., the voltage safety limits) in the OPF problems, thus compromising their practicality. Tackling such inequality constraints is necessary but is challenging for existing machine learning algorithms \cite{eisen2019learning,huang2024unsupervised,gupta2022controlling}. Besides, the incremental control strategy in \cite{yuan2024learning,yuan2024unsupervised} may not generate a steep descent direction toward the varied OPF solutions, as the increment may not accurately follow the gradient of the OPF objectives.

\textit{Contributions:} To overcome the limitations of the existing methods, this paper develops a data-driven real-time algorithm that leverages the available local feedback of DERs to solve time-varying OPF problems. For each time step, the decision of the proposed algorithm is updated based on a term that goes towards the descent direction of the OPF objective, and a term consisting of a learnable function that takes local feedback as input but is trained to capture the system-wide constraint violations, in a way similar to projected gradient descent for constrained optimization. The learnable function, under certain conditions, will stabilize the proposed algorithm to a unique stationary point given an instance of OPF problem. This result motivates us to seek a policy function such that all the stationary points coincide with the (local) optima of a series of OPF problems that change over time, which in turn transfers the original time-varying OPF into a statistical learning problem to optimize the parameters of the learnable function. Remarkably, we show that the time-varying OPF problem can be reformulated as a form that is suitable for conducting statistical learning while satisfying the voltage safety limits with the help of convex approximation to chance constraints \cite{gupta2022controlling,nemirovski2006convex}. The reformulation, with appropriate parametrization of the learnable function, enables us to further develop a stochastic primal-dual learning strategy to directly solve it, bypassing the need for acquiring a training dataset with pre-solved OPF solutions. To summarize, the contributions of this paper are four-fold:
\begin{itemize}
    \item We propose a data-driven real-time algorithm that utilizes the local feedback to continuously track the time-varying OPF solutions. A sufficient condition to attain the unique stationary point is derived to inspire our training strategy for the proposed algorithm (Theorem \ref{thm::1}).
    \item We provide a statistical reformulation for the time-varying OPF problem, and further develop a stochastic primal-dual learning strategy to directly solve the reformulation without resort to ground truth during training.
    \item We advocate a modified deep neural network (DNN) parametrization for the learnable function in the proposed reformulation while satisfying the condition in Theorem \ref{thm::1}, and then establish an OPF solution-tracking bound of the proposed algorithm using universal approximation capability of DNN (Theorem \ref{thm::2}).
    \item A gradient-free learning approach interacting with the nonlinear power flow is developed to compensate for the model linearization error in the algorithm design. 
\end{itemize}

\textit{Outline:} The rest of this paper is organized as follows. Section \ref{sec::OPF} introduces a distribution network model and a time-varying OPF problem formulation. Section \ref{sec::algorithm} proposes a data-driven real-time OPF algorithm and analyzes its equilibria. Section \ref{sec::learning} develops a primal-dual learning strategy to train the algorithm and establishes the tracking error bound. Section \ref{sec::gradfree} describes a gradient-free learning approach. Section \ref{sec::result} reports the numerical results. Section \ref{sec::conclude} concludes the paper.



\section{Model and Problem Formulation}
\label{sec::OPF}

Consider a power distribution network modeled as an undirected tree graph $\mathcal{G}:=\{\mathcal{N}^{+},\mathcal{E}\}$, where $\mathcal{N}^{+}=\{0\}\cup\mathcal{N}$, with $0$ indexing the root node (the substation, also known as the slack bus), and $\mathcal{N}=\{1,\dots,N\}$ indexing the other nodes. The set $\mathcal{E}$ collects all the lines in the tree graph, with $(i,j)\in \mathcal{E}$ representing the line connecting the node pair $i$ and $j$. Assume the power injections (generation minus load) of each node $i\in\mathcal{N}$ can be decomposed as controllable and uncontrollable parts, in which the uncontrollable parts are time-varying. Denote the temporal domain by $\mathcal{T}$ and assume it is discretized by interval $\tau$. For any time instant $t\in\mathcal{T}$, let $p_i^t$ and $q_i^t$ (respectively, $p_{u,i}^t$ and $q_{u,i}^t$) denote the controllable (respectively, uncontrollable) active and reactive power injections of node $i \in\mathcal{N}^+$ at time $t$. Let $v_i^t$ be the squared voltage magnitude at node $i\in\mathcal{N}^+$. For each line $(i,j)\in\mathcal{E}$, let $\ell_{ij}^t$ denote the squared current magnitude, and $P_{ij}^t$ and $Q_{ij}^t$ represent the sending-end active and reactive power, respectively. We denote by $z_{ij}:=r_{ij}+\mathbf{i}x_{ij}$ the series impedance of line $(i,j)\in\mathcal{E}$. 

At time instant $t$, the physical law described by the branch flow model \cite{baran1989optimal,farivar2013branch} is:
\begin{subequations}\label{distflow}
\begin{alignat}{2}
&P_{ij}^t=-p_{u,j}^t-p_j^t+\sum_{k:(j,k)\in\mathcal{E}}P_{jk}^t+r_{ij}\ell_{ij}^t, ~\forall j\in \mathcal{N},\label{distflow::P}\\
&Q_{ij}^t=-q_{u,j}^t-q_j^t+\sum_{k:(j,k)\in\mathcal{E}}Q_{jk}^t+x_{ij}\ell_{ij}^t, ~\forall j\in \mathcal{N},\label{distflow::Q}\\
&v_j^t=v_i^t-2(r_{ij}P_{ij}^t+x_{ij}Q_{ij}^t)+|z_{ij}|^2\ell_{ij}^t, ~\forall (i,j)\in\mathcal{E},\label{distflow::v}\\
&\ell_{ij}^t v_{i}^t=(P_{ij}^t)^2+(Q_{ij}^t)^2, ~\forall (i,j)\in\mathcal{E}.
\end{alignat}
\end{subequations}

Suppose the squared voltage magnitude $v_0^t$ at the slack bus is given and fixed for all time $t$. We define $\boldsymbol{p}^t:=\left[p_1^t,\dots,p_N^t\right]^{\top}$, $\boldsymbol{q}^t:=\left[q_1^t,\dots,q_N^t\right]^{\top}$ as the power control vectors, and let $\boldsymbol{x}^t:=[(\boldsymbol{p}^t)^{\top},(\boldsymbol{q}^t)^{\top}]^{\top}$. The uncontrolled power vectors $\boldsymbol{p}_{u}^t$ and $\boldsymbol{q}_{u}^t$ are defined similarly, and let $\boldsymbol{d}^t:=[(\boldsymbol{p}_{u}^t)^{\top},(\boldsymbol{q}_{u}^t)^{\top}]^{\top}$. Define $\boldsymbol{v}^t:=\left[v^t_1,\dots,v^t_N\right]^{\top}\in\mathbb{R}^N$ as the squared voltage magnitude vector. The implicit function $\boldsymbol{v}^t\left(\boldsymbol{p}^t,\boldsymbol{q}^t,\boldsymbol{d}^t\right)$ 
 is well-defined by the nonlinear branch flow model \eqref{distflow} under normal operating conditions \cite{cwang2018explicit}.

For each node $i \in\mathcal{N}$, assume the controllable power injections are confined by a compact convex set $\mathcal{Y}_i^t$, for instance, a box:
\begin{alignat}{2}\label{region:1}
\mathcal{Y}_i^t=\left\{(p_i,q_i)\mid \underline{p}_i^t\leqslant p_i\leqslant\overline{p}_i^t, \underline{q}_i^t\leqslant q_i\leqslant\overline{q}_i^t\right\}, ~\forall i \in\mathcal{N}.
\end{alignat}

Let $f^t_i\left(p_i^t,q_i^t\right)$ be the time-varying cost function associated with the controllable power injections $\left(p_i^t,q_i^t\right)$ of node $i\in \mathcal{N}$ at time $t$. Without loss of generality, we assume all the cost functions satisfy the following conditions:
\begin{assumption}\label{ass::1}
For any time slot $t\in\mathcal{T}$, the function $f^t_i\left(p_i^t,q_i^t\right), i\in \mathcal{N}$ is continuously differentiable. Moreover, $\boldsymbol{f}^t:=(f_i^t,i\in\mathcal{N}),\forall t\in\mathcal{T}$ is $m$-strongly convex and $\xi$-smooth on $\boldsymbol{\mathcal{Y}}^t:=\prod_{i\in\mathcal{N}}\mathcal{Y}_i^t$, i.e., $\forall \boldsymbol{x},\boldsymbol{y}\in\boldsymbol{\mathcal{Y}}^t,t\in\mathcal{T}$:
\begin{subequations}\label{fproperty}
\begin{alignat}{2}
\left(\nabla \boldsymbol{f}^t(\boldsymbol{x})-\nabla \boldsymbol{f}^t(\boldsymbol{y})\right)^{\top}(\boldsymbol{x}-\boldsymbol{y})&\geqslant m\left\|\boldsymbol{x}-\boldsymbol{y}\right\|_2^2,\label{f:convex}\\
\left\|\nabla \boldsymbol{f}^t(\boldsymbol{x})-\nabla \boldsymbol{f}^t(\boldsymbol{y})\right\|_2&\leqslant \xi \left\|\boldsymbol{x}-\boldsymbol{y}\right\|_2. \label{f:smooth}
\end{alignat}
\end{subequations}
\end{assumption}

Our goal is to solve the time-varying OPF problem:
\begin{subequations}\label{OPF}
\begin{alignat}{2} 
\text{(OPF$^t$) }\min_{\boldsymbol{p}^t,\boldsymbol{q}^t} &\sum_{i\in \mathcal{N}} f_i^t\left(p_i^t,q_i^t\right) \label{OPF::obj}\\
\text{s.t. } \underline{\boldsymbol{v}}&\leqslant \boldsymbol{v}^t(\boldsymbol{p}^t,\boldsymbol{q}^t,\boldsymbol{d}^t)\leqslant\overline{\boldsymbol{v}}\label{OPF::v}\\
&\left(p_i^t,q_i^t\right)\in \mathcal{Y}_i^t, ~\forall i \in\mathcal{N},
\end{alignat}
\end{subequations}
where vectors $\underline{\boldsymbol{v}}$ and $\overline{\boldsymbol{v}}$ are the lower and upper squared voltage safety limits, respectively, and inequality \eqref{OPF::v} is element-wise. Note that the uncontrollable power injections, objectives, and the feasible regions $\mathcal{Y}_i^t,\forall i\in\mathcal{N}$ may change over time. Traditional OPF methods solve the OPF \eqref{OPF} until the algorithm converges and apply the solution to the varied network \cite{peng2016distributed,farivar2013branch}, which may not satisfy the requirements of future power networks with large and fast fluctuations. Real-time OPF methods update the operating setpoint on a faster timescale. Reference \cite{tang2017realtime} designed an L-BFGS-B-based method to track the time-varying OPF solutions, while references \cite{anese2018optimal,tang2018feedback} proposed a primal-dual control strategy using voltage feedback. All of these methods rely on system-wide real-time communications during applications, which poses critical requirements on the communication infrastructure of power systems. Another consequence is that the computation complexity of those methods grows quadratically with the network sizes \cite{zhou2020accelerate}. In the next sections, we will propose a real-time OPF solution strategy, which uses only local measurements.

\section{A Data-Driven Real-Time OPF Algorithm}
\label{sec::algorithm}
The OPF problem \eqref{OPF} is known to be generally nonconvex and hard to be solved by the traditional OPF methods \cite{low2014convex,lavaei2012zeroduality}. To facilitate the design of our algorithm, we introduce a linearization of the branch flow model. 
The terms $\ell_{ij}^t$ in \eqref{distflow::P}-\eqref{distflow::v} are typically much smaller than the other terms and thus can be ignored. This leads to the following linearized power flow model:
\begin{alignat}{2}\label{linearPF}
\boldsymbol{v}^t(\boldsymbol{p}^t,\boldsymbol{q}^t,\boldsymbol{d}^t)=R\boldsymbol{p}^t+X\boldsymbol{q}^t+\boldsymbol{v}_{env}^t,
\end{alignat}
where the matrix $R=\left[R_{ij}\right]_{n\times n}$, and $X=\left[X_{ij}\right]_{n\times n}$, and their entries are calculated by: $R_{ij}:=\sum_{(\zeta,\xi)\in\mathbb{P}_{i\wedge j}}2r_{\zeta\xi}$ and $X_{ij}:=\sum_{(\zeta,\xi)\in\mathbb{P}_{i\wedge j}}2x_{\zeta\xi}$. Here $\mathbb{P}_{i\wedge j}$ denotes 
the common part of the unique paths from nodes $i$ and $j$, respectively, back to the root. The vector $\boldsymbol{v}_{env}^t=v_0^t\boldsymbol{1}+R\boldsymbol{p}_u^t+X\boldsymbol{q}_u^t$ denotes the uncontrollable component of the voltages.

We design our algorithm based on the projected gradient descent methods, while eliminating the system-wide communication by available local measurements. 
For each node $i$ of the power network, let $\hat{v}_i^t$ and $(p_{u,i}^t,q_{u,i}^t)$ be the measurements of the squared voltage magnitude and the uncontrollable power injections, respectively, right before time $t$. That is $\hat{\boldsymbol{v}}^t=\boldsymbol{v}^t(\boldsymbol{p}^{t-1},\boldsymbol{q}^{t-1},\boldsymbol{d}^t)=R\boldsymbol{p}^{t-1}+X\boldsymbol{q}^{t-1}+\boldsymbol{v}_{env}^t$, in which the uncontrolled power injections have changed from $\boldsymbol{d}^{t-1}$ to $\boldsymbol{d}^t$ but the previous controlled power injections are still applied. Our algorithm updates the next operating point as follows for all $i\in \mathcal{N}$:
\begin{subequations}\label{dynamic}
\begin{alignat}{2} 
 p_i^{t}&=\left[p^{t-1}_i-\alpha\left(\nabla_{p_i} f^{t}_i(x_i^{t-1})+u_{\varphi_{p,i}}\left(\hat{v}_i^{t},p^t_{u,i}\right)\right)\right]_{\mathcal{Y}_i^{t}}, \label{dynamic::p} \\
q_i^{t}&=\left[q^{t-1}_i-\alpha\left(\nabla_{q_i} f^{t}_i(x_i^{t-1})+u_{\varphi_{q,i}}\left(\hat{v}_i^{t},q^t_{u,i}\right)\right)\right]_{\mathcal{Y}_i^{t}}, \label{dynamic::q}\\
\boldsymbol{v}^{t}&(\boldsymbol{p}^{t},\boldsymbol{q}^{t},\boldsymbol{d}^t)=R\boldsymbol{p}^{t}+X\boldsymbol{q}^{t}+\boldsymbol{v}_{env}^{t}, \label{dynamic::v}
\end{alignat}
\end{subequations}
where $x_i^t:=[p_i^t,q_i^t]^{\top}$, and $\alpha>0$ is a constant step size. The subscript $[\cdot]_{\mathcal{Y}_i^{t}}$ represents the projections onto the feasible region $\mathcal{Y}_i^t$. Functions $u_{\varphi_{p,i}}\left(\hat{v}_i^{t},p^t_{u,i}\right)$ and $u_{\varphi_{q,i}}\left(\hat{v}_i^{t},q^t_{u,i}\right)$ are learnable functions that take the local measurements as input. Henceforth, we also refer to $u_{\varphi_{p,i}}$, $u_{\varphi_{q,i}}$ as local feedback policies, or parametrizations. Note that both $\nabla_{p_i} f^{t}_i(x_i^{t-1})$ and $u_{\varphi_{p,i}}\left(\hat{v}_i^{t},p^t_{u,i}\right)$ are available locally (similarly for $q_i^t$), and thus the above dynamics is free of communication. Besides, it is easy to verify that the computational complexity of dynamics \eqref{dynamic} grows linearly with the network size, significantly decreasing from the primal-dual type algorithms \cite{anese2018optimal,tang2018feedback,tang2022running} and even those with a hierarchical acceleration \cite{zhou2020accelerate,liang2022hierarchical,liang2024improve}.

We provide a general explanation to the dynamics \eqref{dynamic}. The term $\nabla_{p_i} f^{t}_i(x_i^{t-1})$ is the gradient descent direction of the OPF objective \eqref{OPF::obj}
at time $t$ starting from $x^{t-1}_i$, while the policy $u_{\varphi_{p,i}}$ rectifies the descent direction with consideration to the constraints in \eqref{OPF::v}. Various projected gradient methods can be understood similarly, e.g., the barrier function approach in \cite{gan2016online}, which utilizes a log-barrier function to replace the hard constraints \eqref{OPF::v}, and the primal-dual gradient method in \cite{liang2024improve,anese2018optimal,tang2018feedback,liang2022hierarchical,tang2022running}, which uses dual variables to measure the violations of the constraints \eqref{OPF::v} and determines the descent direction of the primal variables based on the derivatives of both the OPF objective and the dual variables. Our algorithm can also be understood as learning an equilibrium function to control voltages towards the OPF solutions, similarly to references \cite{yuan2024learning,yuan2024unsupervised}. It is worth mentioning that their learnable control strategies either require an OPF solution dataset \cite{yuan2024learning}, or consider no constraint like \eqref{OPF::v} \cite{yuan2024unsupervised}. In the rest of this paper, we will focus on how to learn the functions $u_{\varphi_{p,i}}$ and $u_{\varphi_{q,i}}$ to drive the dynamics \eqref{dynamic} to track the solutions of time-varying OPF problems. 

We first characterize the conditions under which the dynamics \eqref{dynamic} have an equilibrium point. To facilitate our analysis, we rewrite the dynamics \eqref{dynamic} in a vector form. Let $\boldsymbol{u}_{\boldsymbol{\varphi}}:=[\boldsymbol{u}_{\boldsymbol{\varphi}_p}^{\top},\boldsymbol{u}_{\boldsymbol{\varphi}_q}^{\top}]^{\top}$. Defining $\mathcal{A}=[R, X]$, we can rewrite dynamics \eqref{dynamic} compactly as:
\begin{subequations} \label{dynamic:compact}
\begin{alignat}{2}
\boldsymbol{x}^{t}&=\left[\boldsymbol{x}^{t-1}-\alpha\left(\nabla \boldsymbol{f}^t(\boldsymbol{x}^{t-1})+\boldsymbol{u}_{\boldsymbol{\varphi}}(\hat{\boldsymbol{v}}^{t},\boldsymbol{d}^t)\right)\right]_{\boldsymbol{\mathcal{Y}}^t},\label{dycom::x}\\
\boldsymbol{v}^{t}&=\mathcal{A}\boldsymbol{x}^{t}+\boldsymbol{v}_{env}^t.
\end{alignat}
\end{subequations}

In each time slot $t\in\mathcal{T}$, it is convention to suppose the exogenous conditions remain constant, i.e., $\boldsymbol{f}^t$, $\boldsymbol{v}_{env}^t$, $\boldsymbol{\mathcal{Y}}^t$ keep unchanged in time slot $t$, which is necessary for stability analysis \cite{yuan2024learning}. We have the following theorem regarding the existence and uniqueness of the equilibrium of dynamics \eqref{dynamic}:
\begin{theorem}\label{thm::1}
Let Assumption \ref{ass::1} hold. Suppose the function $u_{\varphi_i}:=[u_{\varphi_{p,i}},u_{\varphi_{q,i}}]^{\top}$ satisfies the following conditions, for all $i\in \mathcal{N}$:

\textit{C1).} The function $u_{\varphi_i}(v_i,d_i)$ is separable with respective to variables $v_i$ and $d_i$, i.e., $u_{\varphi_i}(v_i,d_i)=u_{\vartheta_i}(v_i)+u_{\phi_i}(d_i)$;

\textit{C2).} $u_{\vartheta_i}$ is $L_{\vartheta_i}$-Lipschitz continuous, i.e., there exists a constant $L_{\vartheta_i}>0$, such that $\left\|u_{\vartheta_i}(v)-u_{\vartheta_i}(v^{\prime})\right\|_2\leqslant L_{\vartheta_i}\left\|v-v^{\prime}\right\|_2$, $\forall v,v^{\prime}\in \mathbb{R}^{+}$;\\
then the dynamics \eqref{dynamic} has an equilibrium, which must satisfy:\begin{subequations}\label{equilibrium}
\begin{alignat}{2}
\boldsymbol{x}^{\dagger,t}&=\left[\boldsymbol{x}^{\dagger,t}-\alpha\left(\nabla \boldsymbol{f}^t(\boldsymbol{x}^{\dagger,t})+\boldsymbol{u}_{\boldsymbol{\varphi}}(\boldsymbol{v}^{\dagger,t},\boldsymbol{d}^t)\right)\right]_{\mathcal{Y}^t},\label{equilibrium::x}\\
\boldsymbol{v}^{\dagger,t}&=\mathcal{A}\boldsymbol{x}^{\dagger,t}+\boldsymbol{v}_{env}^t.
\end{alignat}
\end{subequations}
Moreover, if $\alpha<\frac{2m}{\xi^2}$ and function $u_{\vartheta_i}$ satisfies:

\textit{C3).} $u_{\vartheta_i}$ has its Lipschitz constant $L_{\boldsymbol{\vartheta}}=\max_{i\in\mathcal{N}} L_{\vartheta_i}<\frac{1-\sqrt{1-2\alpha m+\alpha^2\xi^2}}{\alpha \left\|\mathcal{A}\right\|_2}$;\\
then the dynamics \eqref{dynamic} has a unique equilibrium point.
\end{theorem}

\begin{IEEEproof}
  See Appendix \ref{append::A}.
\end{IEEEproof}

Theorem \ref{thm::1} specifies under which conditions the dynamics \eqref{dynamic} will have a unique equilibrium. Such an equilibrium point depends on the policy function $\boldsymbol{u}_{\boldsymbol{\varphi}}$, may not be optimal, and may violate the voltage constraints \eqref{OPF::v}. Nonetheless, we have an important result directly from Theorem \ref{thm::1}. The equilibrium point $(\boldsymbol{x}^{\dagger,t},\boldsymbol{v}^{\dagger,t})$ is a function of the policy $\boldsymbol{u}_{\boldsymbol{\varphi}}$ that is parametrized by variable $\boldsymbol{\varphi}$. This result implies we can alternatively find a $\boldsymbol{\varphi}$ such that the equilibrium point coincides with the optimal solution of OPF problem \eqref{OPF}. To this end, we define the set $
\boldsymbol{\Phi}=\{\boldsymbol{\varphi}: u_{\varphi_i} \text{ satisfies conditions \textit{C1)-C3)}},~\forall i \in\mathcal{N}\}$ and use the function below to denote the equilibrium \eqref{equilibrium::x} that is uniquely determined by $\boldsymbol{\varphi}$ at time $t$:
\begin{alignat}{2}\label{Ht}
\boldsymbol{x}^{\dagger,t}=H_t(\boldsymbol{\varphi}),\quad\boldsymbol{\varphi}\in\boldsymbol{\Phi}.
\end{alignat}

We have the following lemma characterizing a property of the function $H_t$, which will be useful for our analysis.
\begin{lemma}\label{lemm::1}
The derivative of $H_t$ with respect to $\boldsymbol{u}_{\boldsymbol{\varphi}}$ is bounded $\forall t$, i.e., there exists a constant $L_{h}<\infty$, such that: 
\begin{alignat}{2}
\sup_{t\in\mathcal{T},\boldsymbol{\varphi}\in\boldsymbol{\Phi}}\left\|\frac{\partial H_t(\boldsymbol{\varphi})}{\partial \boldsymbol{u}_{\boldsymbol{\varphi}}}\right\|_2\leqslant L_{h}.
\end{alignat}
\end{lemma}

The proof of Lemma \ref{lemm::1} is merged into the proof of Theorem \ref{thm::2} in Appendix \ref{append::B}.

By using function \eqref{Ht}, the OPF problem \eqref{OPF} at time $t$ becomes the one in which the optimization is over $\boldsymbol{\varphi}$:
\begin{subequations}\label{OPFphi}
\begin{alignat}{2}
\text{(OPF$^t$-$\boldsymbol{\varphi}$) }\min_{\boldsymbol{\varphi}\in\boldsymbol{\Phi}}&\sum_{i\in\mathcal{N}}f^t_i(x_i^{\dagger,t})\label{OPFphi::obj} \\
\text{s.t. } \boldsymbol{x}^{\dagger,t}&=H_t(\boldsymbol{\varphi}),\\
\boldsymbol{v}^{\dagger,t}&(\boldsymbol{\varphi})=\mathcal{A}H_t(\boldsymbol{\varphi})+\boldsymbol{v}_{env}^t,\\
\underline{\boldsymbol{v}}&\leqslant \boldsymbol{v}^{\dagger,t}(\boldsymbol{\varphi})\leqslant \overline{\boldsymbol{v}}, \label{OPFphi::v}
\end{alignat}
\end{subequations}
where $\boldsymbol{v}^{\dagger,t}(\boldsymbol{\varphi})$ is the squared voltage magnitude vector at the equilibrium determined by $\boldsymbol{\varphi}$. Optimizing problem \eqref{OPFphi} may return a solution policy such that the OPF objective attains its optimum while satisfying the voltage constraints. However, such a policy $\boldsymbol{u}_{\boldsymbol{\varphi}}$ may only work for the OPF problem at time $t$, not for the problems at all $t\in\mathcal{T}$. Solving a policy $\boldsymbol{u}_{\boldsymbol{\varphi}}$ for the problem \eqref{OPFphi} at all $t\in\mathcal{T}$ can be challenging, since the optimized policy for the problems may vary with time $t$. If the parametrization of $\boldsymbol{u}_{\boldsymbol{\varphi}}$ is not well selected, i.e., if we utilize a linear parametrization $u_{\varphi_i}=a_iv_i+b_id_i,~\forall i\in\mathcal{N}$, we may not be able to find a solution policy that solves the problem \eqref{OPFphi} at all time $t\in\mathcal{T}$. Another challenge arises from that the problem \eqref{OPFphi} is not well formulated. We may require a specific reformulation to consider the statistical distribution of problem \eqref{OPFphi} in the whole time horizon $\mathcal{T}$. To tackle the challenges, we resort to a class of parametrizations that is near-universal \cite{liang2017why}, which is defined as:
\begin{definition}\label{def::1}
    Let $\boldsymbol{u}_{\boldsymbol{\varphi}^{*,t}}$ be the function of dynamics \eqref{dynamic:compact} that optimizes problem \eqref{OPFphi} at time $t$. For a set of functions $\{\boldsymbol{u}_{\boldsymbol{\varphi}^{*,t}}, t\in\mathcal{T}\}$ for the whole time horizon $\mathcal{T}$, if there exists a parametrization $\boldsymbol{u}_{\boldsymbol{\varphi}}$ with parameter $\boldsymbol{\varphi}\in\boldsymbol{\Phi}$ and a constant $\epsilon>0$ such that
\begin{alignat}{2}
\sup_{t\in\mathcal{T}}\left\|\boldsymbol{u}_{\boldsymbol{\varphi}^{*,t}}-\boldsymbol{u}_{\boldsymbol{\varphi}}\right\|_2\leqslant \epsilon,
\end{alignat}
then $\boldsymbol{u}_{\boldsymbol{\varphi}}$ is called an $\epsilon$-universal parametrization in $\mathcal{T}$.
\end{definition}

Various parametrizations exhibit the universal approximation properties as in Definition \ref{def::1}, e.g., the radial basis function networks \cite{park1991universal} and reproducing kernel Hilbert spaces \cite{sriperumbudur2010relation}. Inspired by the recent success of deep neural networks (DNNs) in functional approximation, we develop our parametrization $\boldsymbol{u}_{\boldsymbol{\varphi}}$ based on DNNs. In the next section, we will introduce our problem reformulation for statistical learning and learning strategy, and elaborate the design of parametrization $\boldsymbol{u}_{\boldsymbol{\varphi}}$ satisfying both Definition \ref{def::1} and conditions \textit{C1)-C3)}. 

\section{Learning Parametrization}
\label{sec::learning}

\subsection{Reformulation via Chance Constraints}
Remember that our goal is to learn a policy $\boldsymbol{u}_{\boldsymbol{\varphi}}$ to drive the designed dynamics \eqref{dynamic} to track the time-varying OPF solutions. To account for the change of the OPF problem \eqref{OPF} and its variant \eqref{OPFphi} over time $t$, we propose a stochastic reformulation that incorporates the OPFs in the whole time horizon $\mathcal{T}$ into a single optimization problem:
\begin{subequations}\label{OPFchance}
\begin{alignat}{2}
\text{(OPF-$\boldsymbol{\varphi}$) }\min_{\boldsymbol{\varphi}\in\boldsymbol{\Phi}}&\mathbb{E}\left[\sum_{i\in\mathcal{N}}f_i(x_i^{\dagger})\right]\label{OPFchance::obj} \\
\text{s.t. } \boldsymbol{x}^{\dagger,t}&=H_t(\boldsymbol{\varphi}), ~t\in\mathcal{T}, \label{OPFchance::x}\\
\boldsymbol{v}^{\dagger,t}&(\boldsymbol{\varphi})=\mathcal{A}H_t(\boldsymbol{\varphi})+\boldsymbol{v}_{env}^t, ~t\in\mathcal{T}, \label{OPFchance::phi}\\
\operatorname{Pr}&\left[v_i^{\dagger}(\boldsymbol{\varphi})\leqslant\underline{v}_i\right]\leqslant \beta,~\forall i\in\mathcal{N}, \label{OPFchance::vlower}\\
\operatorname{Pr}&\left[v_i^{\dagger}(\boldsymbol{\varphi})\geqslant\overline{v}_i\right]\leqslant \beta,~\forall i \in\mathcal{N},\label{OPFchance::vupper}
\end{alignat}
\end{subequations}
where the expectation $\mathbb{E}\left[\cdot\right]$ and the probability $\operatorname{Pr}\left[\cdot\right]$ are with respect to the distribution of $\boldsymbol{f}^t$ and $\boldsymbol{d}^t$ over $t\in\mathcal{T}$, with $\boldsymbol{x}^{\dagger,t},\boldsymbol{v}^{\dagger,t}$ the realization of $\boldsymbol{x}^{\dagger}, \boldsymbol{v}^{\dagger}$ at each $t$. Equations \eqref{OPFchance::x} and \eqref{OPFchance::phi} can be interpreted as a system model that takes the function variable $\boldsymbol{\varphi}$ as input and then outputs an equilibrium point $(\boldsymbol{x}^{\dagger,t},\boldsymbol{v}^{\dagger,t})$ at time $t$. $\operatorname{Pr}[\cdot]$ turns the inequality constraints \eqref{OPFphi::v} into chance constraints. Constant $\beta\in (0,1)$ is a small number that limits the probability of violating the voltage safety limit. An attractive advantage of the formulation \eqref{OPFchance} is that it facilitates the use of prevailing machine learning methods, such as stochastic gradient descent (SGD), to solve for an optimum \cite{eisen2019learning,gupta2022controlling}. Compared to common machine learning methods, our difference lies in the chance constraints \eqref{OPFchance::vlower}-\eqref{OPFchance::vupper}. We rewrite \eqref{OPFchance::vlower}-\eqref{OPFchance::vupper} in the following form to eliminate the probability operator:
\begin{subequations}\label{OPFchanceI}
\begin{alignat}{2}
\mathbb{E}&\left[\mathds{1}\left(\underline{v}_i-v_i^{\dagger}(\boldsymbol{\varphi})\right)\right]\leqslant \beta,~\forall i\in\mathcal{N}, \label{OPFchanceI::vlower}\\
\mathbb{E}&\left[\mathds{1}\left(v_i^{\dagger}(\boldsymbol{\varphi})-\overline{v}_i\right)\right]\leqslant \beta,~\forall i \in\mathcal{N},\label{OPFchanceI::vupper}
\end{alignat}
\end{subequations}
where function $\mathds{1}(x)$ is an indicator that equals $1$ if $x\geqslant 0$ and $0$ otherwise.
However, the expectation constraints \eqref{OPFchanceI} are still difficult to address, since the indicator function is nonconvex and non-differentiable, thus obstructing the calculation of gradients. 

\begin{figure}
    \centering
    \includegraphics[width=1.0\columnwidth]{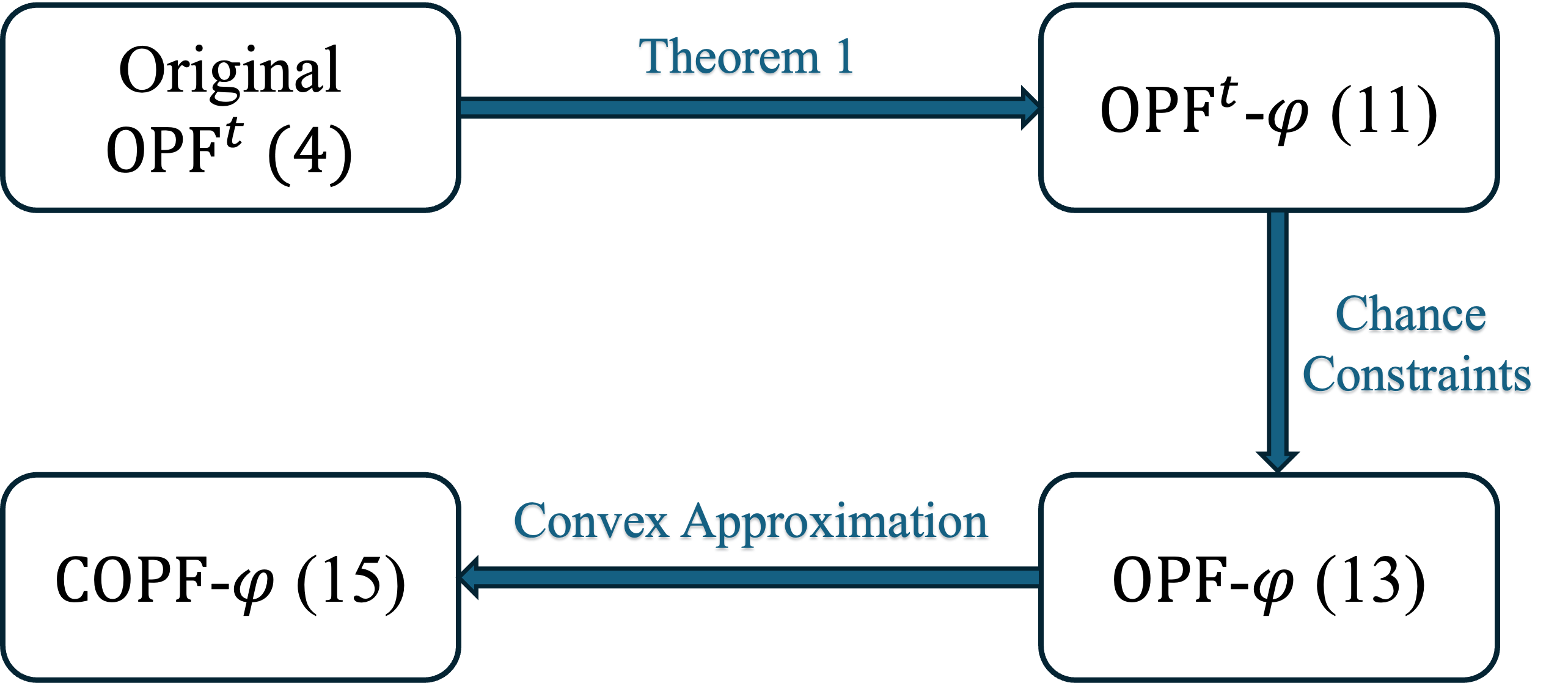}
    \caption{Illustration of OPF problem transformations for learning.}
    \label{fig:transformation}
\end{figure}

We briefly review an existing method for convex approximation of chance constraints \cite{gupta2022controlling,nemirovski2006convex}. Consider some function $g(x)$ of random variable $x$ and chance constraint $\operatorname{Pr}\left[g(x)\geqslant 0\right]=\mathbb{E}\left[\mathds{1}\left(g(x)\right)\right]\leqslant \beta$. It is easy to verify $\mathds{1}\left(g(x)\right)\leqslant\left[1+g(x)/\lambda\right]_{+}$, where $[\cdot]_{+}$ denotes the projection onto the nonnegative orthant, for any $g(x)$ and $\lambda>0$. Therefore, if $\mathbb{E}\left[\left[1+g(x)/\lambda\right]_{+}\right]\leqslant \beta$ for $\beta\in(0,1)$, then the chance constraint $\operatorname{Pr}\left[g(x)\geqslant0\right]\leqslant\beta$ holds, too. Multiplying both sides of $\mathbb{E}\left[\left[1+g(x)/\lambda\right]_{+}\right]\leqslant \beta$ by $\lambda$ and using $\lambda>0$, we have a convex surrogate $\mathbb{E}\left[\left[\lambda+g(x)\right]_{+}\right]\leqslant\beta \lambda$. This result can be further extended to $\lambda\in\mathbb{R}$; see \cite{gupta2022controlling,nemirovski2006convex} and references therein. Replacing $g(x)$ with $\underline{v}_i-v_i^{\dagger}(\boldsymbol{\varphi})$ and $v_i^{\dagger}(\boldsymbol{\varphi})-\overline{v}_i$ respectively and collecting the corresponding constraints over $i\in\mathcal{N}$ in a vector form, we have the following convex approximation to problem \eqref{OPFchance}:
\begin{subequations}\label{OPFchanceII}
\begin{alignat}{2}
\text{(COPF-$\boldsymbol{\varphi}$)}\min_{\boldsymbol{\varphi}\in\boldsymbol{\Phi},\underline{\boldsymbol{\lambda}},\overline{\boldsymbol{\lambda}}}&\mathbb{E}\left[\sum_{i\in\mathcal{N}}f_i(x_i^{\dagger})\right]\label{OPFchanceII::obj} \\
\text{s.t. } \boldsymbol{x}^{\dagger,t}&=H_t(\boldsymbol{\varphi}), ~t\in\mathcal{T}, \label{OPFchanceII::x}\\
\boldsymbol{v}^{\dagger,t}&(\boldsymbol{\varphi})=\mathcal{A}H_t(\boldsymbol{\varphi})+\boldsymbol{v}_{env}^t,~t\in\mathcal{T},\label{OPFchanceII::phi}\\
\mathbb{E}&\left[\left[\underline{\boldsymbol{\lambda}}+\underline{\boldsymbol{v}}-\boldsymbol{v}^{\dagger}(\boldsymbol{\varphi})\right]_{+}\right]\leqslant \beta \underline{\boldsymbol{\lambda}}, \label{OPFchanceII::vlower}\\
\mathbb{E}&\left[\left[\overline{\boldsymbol{\lambda}} +\boldsymbol{v}^{\dagger}(\boldsymbol{\varphi})-\overline{\boldsymbol{v}}\right]_{+}\right]\leqslant \beta \overline{\boldsymbol{\lambda}},\label{OPFchanceII::vupper}
\end{alignat}
\end{subequations}
where $\underline{\boldsymbol{\lambda}}=[\underline{\lambda}_1,\dots,\underline{\lambda}_N]^{\top}\in\mathbb{R}^N$ and $\overline{\boldsymbol{\lambda}}=[\overline{\lambda}_1,\dots,\overline{\lambda}_N]^{\top}\in\mathbb{R}^N$ are auxiliary variables. The formulation \eqref{OPFchanceII} permits the SGD-based solution strategy to statistically learn a solution $\boldsymbol{\varphi}$ with respect to the distribution over $\mathcal{T}$. The problem transformations from the original OPF \eqref{OPF} to the chance-constrained formulation \eqref{OPFchanceII} are summarized in Figure \ref{fig:transformation}. Note that these may not be equivalent transformations. However, if the learned parametrization $\boldsymbol{u}_{\boldsymbol{\varphi}}$ satisfies Definition \ref{def::1}, the OPF solution-tracking error of dynamics \eqref{dynamic} can be bounded as we will show shortly.

\subsection{Primal-Dual Learning}
Traditional machine learning objectives are typically optimized in an unconstrained manner via SGD-based algorithms, while problem \eqref{OPFchanceII} involves both equality and inequality constraints and thus cannot be solved by conventional machine learning algorithms. As mentioned before, we treat the equations \eqref{OPFchanceII::x}-\eqref{OPFchanceII::phi} as a system model that admits a function variable $\boldsymbol{\varphi}$. Then the equilibrium point $(\boldsymbol{x}^{\dagger,t},\boldsymbol{v}^{\dagger,t})$ will be automatically solved by the system model. Such a model is realizable using a digital twin of the power network (e.g., a software simulator such as Matpower or OpenDSS) or a hardware simulator. Moreover, the communication needed for solving problem \eqref{OPFchanceII} can be done in the digital twin or hardware simulator in a cost-effective manner, while the application of dynamics \eqref{dynamic} to the real power network for real-time optimization needs no more communication. The major challenge remaining is the chance constraints \eqref{OPFchanceII::vlower}-\eqref{OPFchanceII::vupper}, which can be addressed by formulating the Lagrangian function of problem \eqref{OPFchanceII} and applying the stochastic primal-dual updates \cite{eisen2019learning} to approach a solution. Let $\underline{\boldsymbol{\mu}}$ and $\overline{\boldsymbol{\mu}}$ be the dual variables associated with the inequality constraints \eqref{OPFchanceII::vlower} and \eqref{OPFchanceII::vupper}, respectively. Consider the Lagrangian function of problem \eqref{OPFchanceII}:
\begin{alignat}{2}\label{lagrangian}
\mathcal{L}(\boldsymbol{\varphi},\boldsymbol{\lambda},\boldsymbol{\mu})=\underline{\boldsymbol{\mu}}^{\top}\left(\mathbb{E}\left[[\underline{\boldsymbol{\lambda}}+\underline{\boldsymbol{v}}-\boldsymbol{v}^{\dagger}(\boldsymbol{\varphi})]_{+}\right]-\beta\underline{\boldsymbol{\lambda}}\right)+\nonumber\\
\overline{\boldsymbol{\mu}}^{\top}\left(\mathbb{E}\left[[\overline{\boldsymbol{\lambda}} +\boldsymbol{v}^{\dagger}(\boldsymbol{\varphi})-\overline{\boldsymbol{v}}]_{+}\right]-\beta\overline{\boldsymbol{\lambda}}\right)+\mathbb{E}\Big[\sum_{i\in\mathcal{N}}f_i(x_i^{\dagger})\Big],
\end{alignat}
where $\boldsymbol{\lambda}:=[\underline{\boldsymbol{\lambda}}^{\top},\overline{\boldsymbol{\lambda}}^{\top}]^{\top}$ and $\boldsymbol{\mu}:=[\underline{\boldsymbol{\mu}}^{\top},\overline{\boldsymbol{\mu}}^{\top}]^{\top}$. Function \eqref{lagrangian} replaces the inequality constraints \eqref{OPFchanceII::vlower} and \eqref{OPFchanceII::vupper} by penalizing their violations. Solving for a local optimum of problem \eqref{OPFchanceII} is equivalent to finding a saddle point of the unconstrained max-min problem:
\begin{alignat}{2}\label{minmax}
\max_{\boldsymbol{\mu}\geqslant0}\min_{\boldsymbol{\varphi}\in\boldsymbol{\Phi},\boldsymbol{\lambda}}\mathcal{L}(\boldsymbol{\varphi},\boldsymbol{\lambda},\boldsymbol{\mu}).
\end{alignat}

In particular, a stochastic primal-dual gradient algorithm \cite{eisen2019learning} to approach a saddle point of \eqref{minmax} is:
\begin{subequations}\label{primaldual}
\begin{alignat}{2}
    &\boldsymbol{\varphi}^{k}=\left[\boldsymbol{\varphi}^{k-1}-\sigma_{\varphi}\nabla_{\boldsymbol{\varphi}}\mathcal{L}(\boldsymbol{\varphi}^{k-1},\boldsymbol{\lambda}^{k-1},\boldsymbol{\mu}^{k-1})\right]_{\boldsymbol{\Phi}}, \label{primal::phi}\\
    &\underline{\boldsymbol{\lambda}}^{k}=\underline{\boldsymbol{\lambda}}^{k-1}-\sigma_{\lambda}\nabla_{\underline{\boldsymbol{\lambda}}}\mathcal{L}(\boldsymbol{\varphi}^{k-1},\boldsymbol{\lambda}^{k-1},\boldsymbol{\mu}^{k-1}),\label{primal::ulambdau}\\
    &\overline{\boldsymbol{\lambda}}^{k}=\overline{\boldsymbol{\lambda}}^{k-1}-\sigma_{\lambda}\nabla_{\overline{\boldsymbol{\lambda}}}\mathcal{L}(\boldsymbol{\varphi}^{k-1},\boldsymbol{\lambda}^{k-1},\boldsymbol{\mu}^{k-1}),\label{primal::olambdau}\\
    &\underline{\boldsymbol{\mu}}^{k}=\Big[\underline{\boldsymbol{\mu}}^{k-1}+\sigma_{\mu}\mathbb{E}\big[[\underline{\boldsymbol{\lambda}}^{k}+\underline{\boldsymbol{v}}-\boldsymbol{v}^{\dagger}(\boldsymbol{\varphi}^{k})]_{+}-\beta\underline{\boldsymbol{\lambda}}^{k}\big]\Big]_{+},\\
    &\overline{\boldsymbol{\mu}}^{k}=\Big[\overline{\boldsymbol{\mu}}^{k-1}+\sigma_{\mu}\mathbb{E}\big[[\overline{\boldsymbol{\lambda}}^{k}+\boldsymbol{v}^{\dagger}(\boldsymbol{\varphi}^{k})-\overline{\boldsymbol{v}}]_{+}-\beta\overline{\boldsymbol{\lambda}}^{k}\big]\Big]_{+},
\end{alignat}   
\end{subequations}
where $(\sigma_{\varphi},\sigma_{\lambda})$ and $\sigma_{\mu}$ are positive step sizes for the primal and dual updates, respectively. 

We next compute the needed gradients for the primal updates \eqref{primal::phi}-\eqref{primal::olambdau}. Using the rule in stochastic gradient descent, we unfold the expectation operator in \eqref{lagrangian} as the average of batched samples $\{(\boldsymbol{x}^{\dagger,s},\boldsymbol{v}^{\dagger,s})\}_{s=1}^{S}$ from the distribution over $\mathcal{T}$. Since $\boldsymbol{\varphi}:=[\boldsymbol{\varphi}_{p}^{\top},\boldsymbol{\varphi}_{q}^{\top}]^{\top}$, we only elaborate the gradient associated with parameter $\varphi_{p,i}$ at node $i$. The gradient related to $\boldsymbol{\varphi}_{q}$ can be calculated similarly. Using the chain rule of derivatives, we have: 
\begin{alignat}{2}\label{grad::ellphi}
&\frac{\partial \mathcal{L}}{\partial \varphi_{p,i}}=\frac{1}{S}\sum_{s=1}^{S}\left(-\sum_{j=1}^{N}\underline{\mu}_j\frac{\partial v_j^{\dagger,s}}{\partial p_i^{\dagger,s}}\frac{\partial p_i^{\dagger,s}}{\partial \varphi_{p,i}}\mathds{1}\left(\underline{\lambda}_j+\underline{v}_j-v_j^{\dagger,s}\right)+\right.\nonumber\\
&\qquad\left.\sum_{j=1}^{N}\overline{\mu}_j\frac{\partial v_j^{\dagger,s}}{\partial p_i^{\dagger,s}}\frac{\partial p_i^{\dagger,s}}{\partial \varphi_{p,i}}\mathds{1}\left(\overline{\lambda}_j+v_j^{\dagger,s}-\overline{v}_j\right)+\frac{\partial f_i^s}{\partial p_{i}^{\dagger,s}}\frac{\partial p_{i}^{\dagger,s}}{\partial \varphi_{p,i}}\right)\nonumber\\
&\quad=\frac{1}{S}\sum_{s=1}^{S}\left[R_{ij}\left(-\sum_{j=1}^{N}\underline{\mu}_j \mathds{1}\left(\underline{\lambda}_j+\underline{v}_j-v_j^{\dagger,s}\right)+\right.\right.\nonumber\\
&\quad\left.\left.\quad\sum_{j=1}^{N}\overline{\mu}_j \mathds{1}\left(\overline{\lambda}_j+v_j^{\dagger,s}-\overline{v}_j\right)\right)+\frac{\partial f_i^s}{\partial p_i^{\dagger,s}}\right]\cdot\frac{\partial p_i^{\dagger,s}}{\partial \varphi_{p,i}},
\end{alignat}
where we calculate $\partial v_j^{\dagger,s}/\partial p_i^{\dagger,s}=R_{ij}$ in the second equality using the linearized model \eqref{linearPF}. The error introduced by the linearization was shown to be small \cite{anese2018optimal} and can be bounded with respect to a local optimum of the actual nonconvex OPF problem \cite{liang2024improve,liang2022hierarchical}. The term $\partial f_i^s/\partial p_i^{\dagger,s}$ depends on the specific OPF objective and is easy to compute. The remaining term $\partial p_i^{\dagger,s}/\partial \varphi_{p,i}$, by taking derivatives on both sides of \eqref{equilibrium::x}, can be calculated as:
\begin{alignat}{2}
\frac{\partial p_i^{\dagger,s}}{\partial \varphi_{p,i}}=\begin{cases}-\frac{1}{\nabla_{p_ip_i}^2 f^s_i}\frac{\partial u_{\varphi_{p,i}}(v_i^{\dagger,s},p^s_{u,i})}{\partial \varphi_{p,i}}, & g_i^s(p_i^{\dagger,s})\in\mathcal{Y}_i^s,  \\ 0, & \text{otherwise} ,\end{cases}
\end{alignat}
where $g_i(p_i^{\dagger})=p_i^{\dagger}-\alpha(\nabla_{p_i}f_i(p_i^{\dagger},q_i^{\dagger})+u_{\varphi_{p,i}}(v_i^{\dagger},p_{u,i}))$. The term $\partial u_{\varphi_{p,i}}(v_i^{\dagger,s},p^s_{u,i})/ \partial \varphi_{p,i}$ is the partial derivative of the output with respect to the parameters of the learnable function $u_{\varphi}$, which can be readily calculated using the backpropagation of the learning process. The gradients $\partial \mathcal{L}/\partial \boldsymbol{\lambda}$ in \eqref{primal::ulambdau}-\eqref{primal::olambdau} can be calculated as:
\begin{subequations}\label{grad::elllambda}
\begin{alignat}{2}
\frac{\partial \mathcal{L}}{\partial \underline{\lambda}_i}&=\underline{\mu}_i\cdot\frac{1}{S}\sum_{s=1}^{S}\left(\mathds{1}\left(\underline{\lambda}_i+\underline{v}_i-v_i^{\dagger,s}\right)-\beta\right),\\
\frac{\partial\mathcal{L}}{\partial \overline{\lambda}_i}&=\overline{\mu}_i\cdot\frac{1}{S}\sum_{s=1}^{S}\left(\mathds{1}\left(\overline{\lambda}_i+v_i^{\dagger,s}-\overline{v}_i\right)-\beta\right).
\end{alignat}
\end{subequations}

\begin{figure}
    \centering
    \includegraphics[width=1.0\columnwidth]{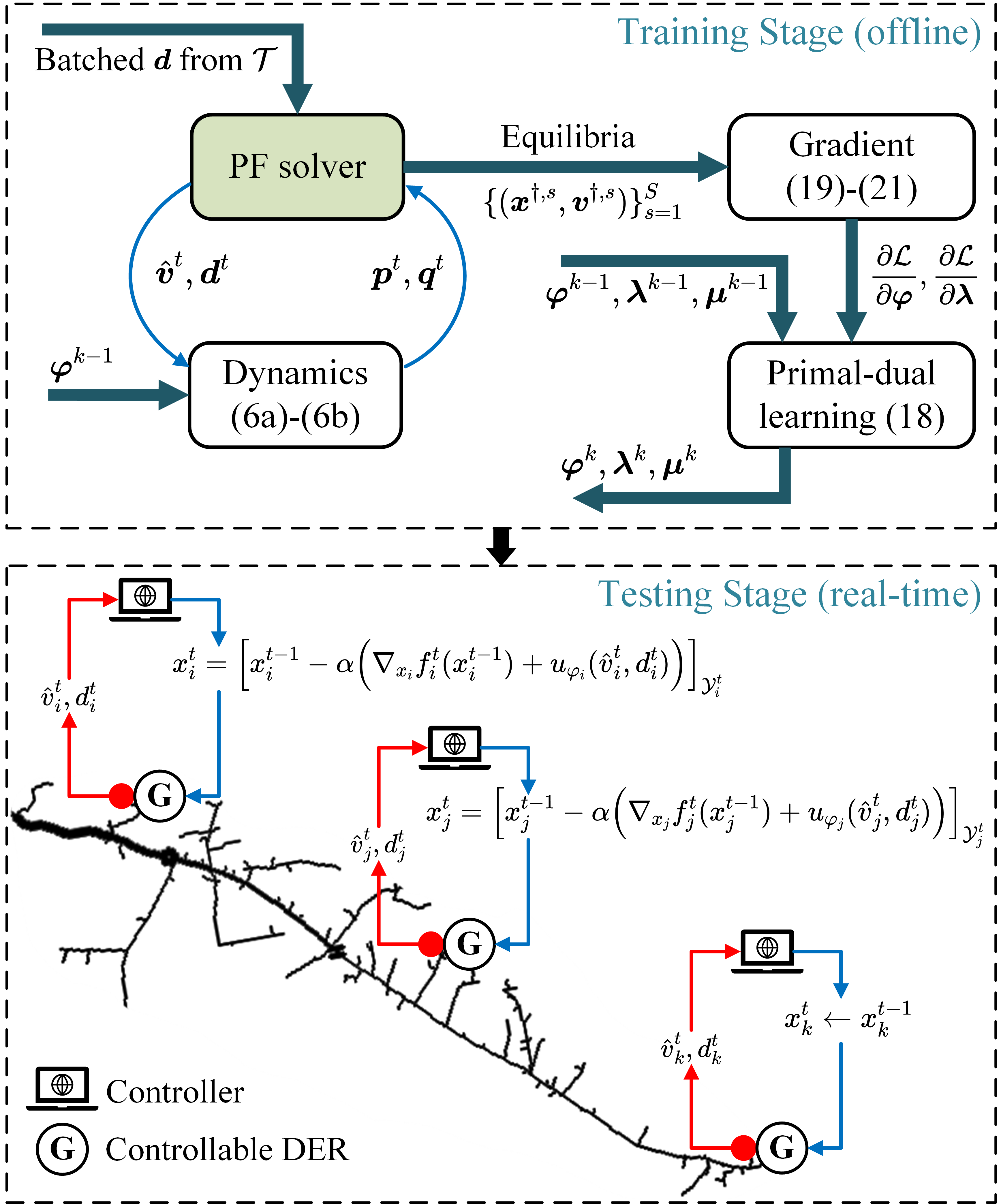}
    \caption{Framework of training and testing (operation) stages.}
    \label{fig:learning}
\end{figure}

We refer to the stochastic primal-dual updates \eqref{primaldual} with the gradients \eqref{grad::ellphi}-\eqref{grad::elllambda} over batches of samples from the distribution over $\mathcal{T}$ as the training (learning) stage, and the application of dynamics \eqref{dynamic} with the learned solution $\boldsymbol{\varphi}$ to real power networks as the testing (operation) stage. The overall framework of the proposed method is summarized in Figure \ref{fig:learning}. We highlight that we need not know the exact distribution over $\mathcal{T}$, as acquiring the power injections from future time $t$ is unrealistic. Instead, we can use the historical samples or the predictions of uncontrollable injections as the training set for the primal-dual learning.

\subsection{DNN Parametrization and Tracking Performance}
We have so far discussed the primal-dual learning technique to find a solution $\boldsymbol{\varphi}$. In this part, we further develop a parametrization $\boldsymbol{u}_{\boldsymbol{\varphi}}$ that carries the near-universal approximation properties as in Definition \ref{def::1}, while satisfying the conditions \textit{C1)-C3)}. In particular, we start our design from the popular deep neural networks (DNNs) that have exhibited great performance in function approximation. We then establish the error bounds for the OPF solution-tracking performance of the proposed dynamics \eqref{dynamic} by virtue of the $\epsilon$-universal property in Definition \ref{def::1}. We only elaborate the design of $u_{\varphi_{p,i}}$ at node $i\in\mathcal{N}$ and that for $u_{\varphi_{q,i}}$ is similar. The parametrization of $\varphi_{p,i}$ is modeled by a single-input single-output fully connected neural network that carries $p_{u,i}$ plus a learnable monotone hyper-connection $k_{p,i}v_i$:\begin{subequations}\label{NN}
  \begin{alignat}{2}
   h_0&=p_{u,i},\\
   h_{\ell}&=\sigma_{\ell}(\mathbf{W}_{\ell}h_{\ell-1}+\mathbf{b}_{\ell}), ~\forall \ell=1,\dots,L,\\
   u_{\varphi_{p,i}}&=\mathbf{W}_{L+1} h_{L}+\mathbf{b}_{L+1}+k_{p,i}v_i,
\end{alignat}  
\end{subequations}
where $\{\mathbf{W}_{\ell}\in\mathbb{R}^{n_{\ell}\times n_{\ell-1}}\}_{\ell=1}^{L+1}$ and $\{\mathbf{b}_{\ell}\in\mathbb{R}^{n_{\ell}}\}_{\ell=1}^{L+1}$ are the learnable parameters for input $p_{u,i}$ with $n_{\ell}$ denoting the dimension of layer $\ell$ and $L$ denoting the number of hidden layers, and $k_{p,i}$ is a learnable scalar parameter for input $v_i$. All these parameters are contained in $\varphi_{p,i}$. Activation function $\sigma_{\ell}$ is the ReLU function for all hidden layers. The conditions \textit{C1)-C3)} are enforced by directly feeding $v_i$ through a linear term $k_{p,i}v_i$ to the output layer.

\begin{figure}
    \centering
    \includegraphics[width=0.8\columnwidth]{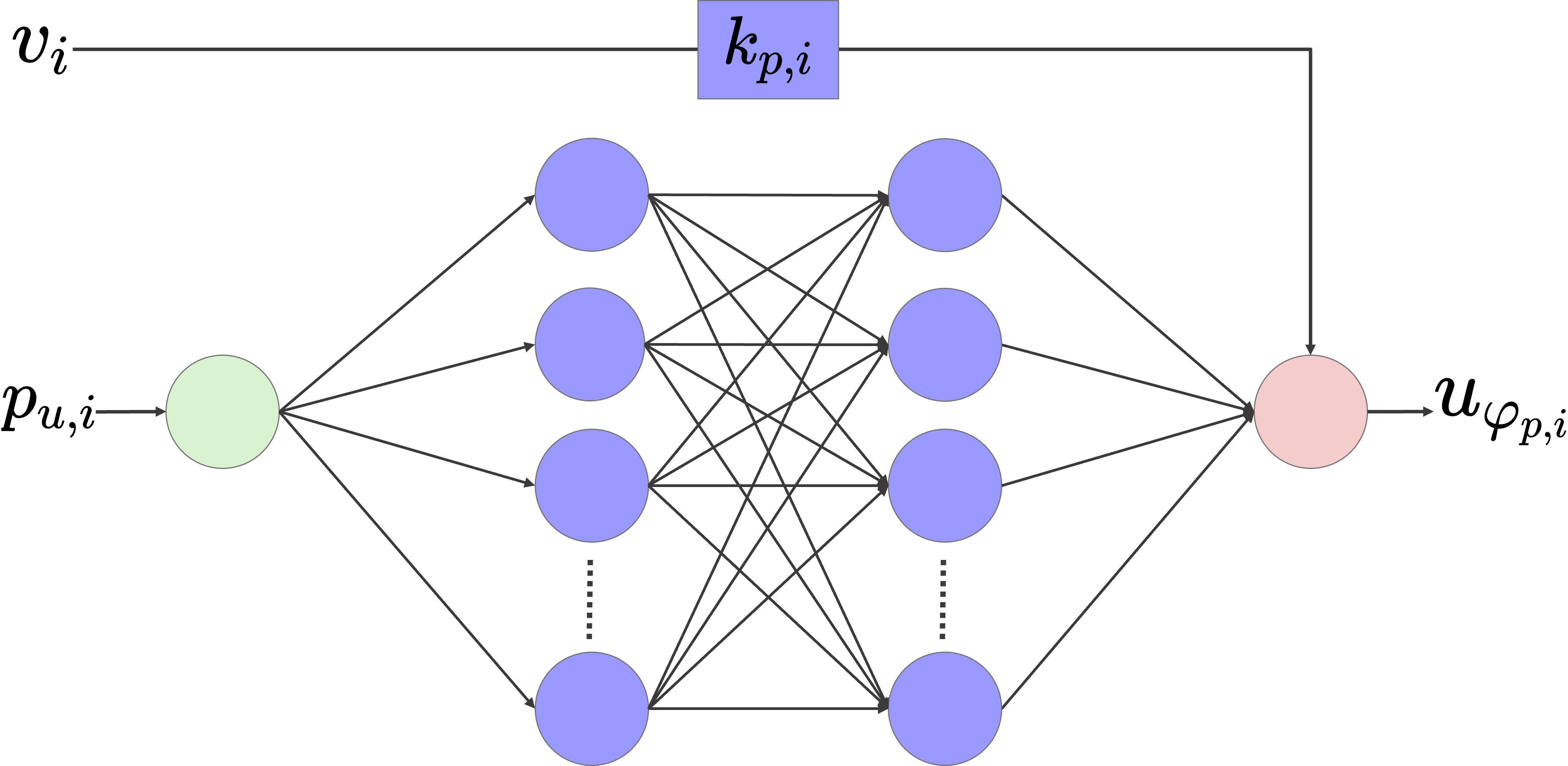}
    \caption{Architecture of the proposed neural network.}
    \label{fig:neuralnetwork}
\end{figure}

We demonstrate the designed parametrization \eqref{NN} in Figure \ref{fig:neuralnetwork}. Indeed, the approximation capability of the deep neural network \eqref{NN} grows with the number of layers $L$ and the layer sizes $n_{\ell}$. The needed number of layers and layer sizes to attain the desired approximation accuracy $\epsilon$ in Definition \ref{def::1} were studied in \cite{liang2017why,yarotsky2017errorbound}, which control the tradeoff between the approximation accuracy and the computational efficiency. 

With this parametrization, we next characterize the OPF solution-tracking error bound of the dynamics \eqref{dynamic}. We define:\begin{alignat}{2}
\rho(\alpha)=&\big[1+\alpha^2(\xi^2+L_{\boldsymbol{\vartheta}}^2\left\|\mathcal{A}\right\|^2_2+2\xi L_{\boldsymbol{\vartheta}}\left\|\mathcal{A}\right\|_2)-2\alpha m\big]^\frac{1}{2}.
\end{alignat}

Let $\boldsymbol{x}^{*,t}$ denote the unique optimizer of OPF problem \eqref{OPF} with the linearized model at time $t$. We denote:
\begin{alignat}{2}\label{xstarchange}
    \gamma=\sup_{t\in\mathcal{T}/\{0\}}\left\|\boldsymbol{x}^{*,t}-\boldsymbol{x}^{*,t-1}\right\|_2.
\end{alignat}



The following theorem establishes the tracking error bound of dynamics \eqref{dynamic} with local feedback $(\hat{v}^t_i,d^t_i),\forall i\in \mathcal{N},t\in\mathcal{T}$. %
\begin{theorem} \label{thm::2}
    Suppose the parametrization $\boldsymbol{u}_{\boldsymbol{\varphi}^{*}}$ returned by \eqref{primaldual} is $\epsilon$-universal as defined in Definition \ref{def::1}. Let Assumption \ref{ass::1} and conditions \textit{C1)}-\textit{C3)} hold. If the function $u_{\vartheta_i}(\cdot)$ is non-decreasing in $v_i$, and the step size $\alpha>0$ is chosen such that:
    $$\rho(\alpha)<1,$$
    then the sequence $\{\boldsymbol{x}^t\}_{t\in\mathcal{T}}$ generated by the data-driven dynamics \eqref{dynamic} converges exponentially to $\boldsymbol{x}^{*,t}$ up to the asymptotic bound given by:
    \begin{alignat}{2}\label{bound}
        \limsup_{t\rightarrow \infty}\left\|\boldsymbol{x}^{t}-\boldsymbol{x}^{*,t}\right\|_2=\frac{\rho(\alpha)\gamma+(1+\rho(\alpha))L_{h}\epsilon}{1-\rho(\alpha)}.
    \end{alignat}
\end{theorem}
\begin{IEEEproof}
    See Appendix \ref{append::B}.
\end{IEEEproof}

We now analyze the tracking error bound in \eqref{bound}. The first term is proportional to $\gamma$, the maximum rate of change of the optimizer $\boldsymbol{x}^{*,t}$ over time, which is common in time-vary optimization algorithms \cite{tang2018feedback,tang2022running,anese2018optimal}. The second term is proportional to $L_h$ and the parameter $\epsilon$. This term is characterized by the approximation capability of the parametrization $\boldsymbol{u}_{\boldsymbol{\varphi}}$ for the primal-dual learning process. 
Generally, utilizing deeper neural networks can reduce the second term, while the increased computation burden of DNNs necessitates longer solution time for dynamics \eqref{dynamic}, which potentially results in a larger $\gamma$ and thus increases the first term. Such a tradeoff in the design of parametrization $\boldsymbol{u}_{\boldsymbol{\varphi}}$ can be tuned according to the practical operation requirements of power systems.

\section{Gradient-Free Learning with Nonlinear Model}
\label{sec::gradfree}
In this section, we implement both the training and testing with the nonlinear power flow model \eqref{distflow}. Still, the update of the operating point follows the rules in \eqref{dynamic::p}-\eqref{dynamic::q}, while we utilize the voltage measurements from the nonlinear power flow model \eqref{distflow} rather than the simplified linearized model \eqref{dynamic::v}. 
In fact, utilizing the feedback from the realistic nonlinear power flow can partly compensate for the modeling error underlying the gradient calculation, e.g., the error introduced by $R_{ij}$ in \eqref{grad::ellphi} based on the linearized power flow \cite{anese2018optimal}. However, even with such feedback, calculating the gradients still requires the exact knowledge of the feeder and thus tends to be inaccurate. Noticing that estimating the gradients is only needed in the learning stage, we develop a gradient-free learning approach, which is more practical than computing the gradients in the nonlinear power flow model. The gradient-free approach is especially useful if the power network is large and equipped with devices that have complicated models, e.g., transformers, regulators, capacitors. The existing methods may fail to return exact gradient values for such cases \cite{gan2016online}. Further, in unbalanced three-phase networks, the mutual coupling of phases obstructs the computation of gradients. To circumvent the difficulties in computing the gradient $\partial \boldsymbol{v}/\partial \boldsymbol{x}$ in \eqref{grad::ellphi} from the nonlinear power flow model, we leverage the zero-order optimization technique to construct a finite-difference gradient estimator, by querying the feeder with random perturbations. Specifically, we construct the $2n$-point gradient estimator around the equilibrium $(\boldsymbol{x}^{\dagger},\boldsymbol{v}^{\dagger})$ as:
\begin{alignat}{2}\label{grad::free}
\hat{\nabla}_{\boldsymbol{x}}\boldsymbol{v}(\boldsymbol{x}^{\dagger})= \sum_{k=1}^{n}\frac{\boldsymbol{v}(\boldsymbol{x}^{\dagger}+\varepsilon e_k)-\boldsymbol{v}(\boldsymbol{x}^{\dagger}-\varepsilon e_k)}{2\varepsilon}e_k, 
\end{alignat}
where $\varepsilon$ is a positive number, $e_k$ is the basis vector with all the entries being $0$ except the $k$-th being $1$. $n$ is the dimension of $\boldsymbol{x}$. Indeed, the estimated gradient $\hat{\nabla}_{\boldsymbol{x}}\boldsymbol{v}(\boldsymbol{x}^{\dagger})$ in \eqref{grad::free} can be arbitrarily close to the gradient of the nonlinear power flow model \cite{tang2021distributed}. Since the primal-dual learning stage is conducted offline, the increased computation introduced by the $2n$-point gradient estimator \eqref{grad::free} does not impact the real-time application of the dynamics \eqref{dynamic}. Therefore, with \eqref{grad::free}, we can learn a more reliable model without sacrificing its applicability and responsiveness.

\section{Numerical Results}
\label{sec::result}
We build a single-phase version of IEEE 37-node test feeder, by averaging the line parameters and default loads over three phases of the original data on IEEE PES website (https://cmte.ieee.org/pes-testfeeders/resources/). Denote the single-phase default load at each node $i$ by $d_i^{\text{def}}$. We select $13$ nodes to add controllable generation resources and deploy net loads (i.e., loads minus uncontrollable renewable generations) $\kappa^t_id_i^{\text{def}}$ varying over time $t$; the remaining nodes just have fixed uncontrollable loads $d_i^{\text{def}}$ (which may be zero) all the time. The time-varying factor $\kappa_i^t=\kappa^t_{\text{CA}}+\kappa^t_{\text{rnd},i}$ is designed as follows. First, we interpolate the $5$-minute \emph{net} demand data of California Independent System Operator (CAISO) \cite{california} from 16:00 to 24:00 on 09/09/2024 into $6$-second resolution and divide the data by its maximum value to obtain time-varying ratio $\kappa_{\text{CA}}^t\leqslant 1$. The \emph{net} demand equals the system demand minus the wind and solar generations. Under this setting, the solar generation starts to decline at 16:00, and hence gradually more demand is supplied by the controllable generations. Second, we generate a Gaussian random disturbance $\kappa^t_{\text{rnd},i}=\frac{1}{\sqrt{d_i^{\text{def}}}}\mathbf{N}(1,0.1)$ independently at each node $i$ and time $t$. In this way, we obtain time-varying net loads that incorporate both the slower trend and the faster random disturbances in user demand and renewable generation.

For our experiments, we fix the voltage magnitude of the root (slack) node at $1$ per unit (p.u.), and set the safe voltage limits at $\underline{v}=0.95^2$ p.u. and $\overline{v}=1.05^2$ p.u., respectively. Each controllable power source has the active power capacity $\overline{p}_i^t=500$ kVA and reactive power capacity $\overline{q}_i^t=300$ kVAR, while the lower limits are $\underline{p}_i^t=0,\underline{q}_i^t=0$. The objective function in the time-varying OPF problem \eqref{OPF} is $f^t_i(p_i^t,q_i^t)=(p_i^t-\underline{p}_i^t)^2+(q_i^t-\underline{q}_i^t)^2$, which measures the cost of generation. MATPOWER is used to simulate the nonlinear power flow.

For the proposed algorithm, we build on Pytorch the DNNs parametrized by $\varphi_{p,i}$ and $\varphi_{q,i}$, which have three layers of size 64. The real-time OPF algorithm updates the operating point every $6$ seconds with a step size $\alpha=0.48$. We utilize the net demand data from CAISO on 03/09, 05/09, 08/09/2024, which are processed in the same way as those on 09/09/2024, as the training dataset, and use the data on 09/09/2024 for testing. We apply the primal-dual updates \eqref{primaldual} to train the DNNs. The Adam optimizer with a learning rate $\sigma_{\varphi}=1\times10^{-3}$ is used for the training of the DNNs. Auxiliary variables for chance constraints are empirically fixed as $\overline{\lambda}_i=\underline{\lambda}_i=5\times10^{-4},\forall i\in\mathcal{N}$ to ease the training. The dual variables are updated using SGD with a learning rate $\sigma_{\mu}=100$. We use a minibatch size of $32$ and train the DNNs for $50$ epochs. All the experiments are run in Python 3.8 programs on a Macbook Pro with $8$-core M1 Pro CPU, 16GB RAM, and MacOS.

\subsection{Gradient-Based Learning}
\label{subsec::based}
\begin{figure}
    \centering
    \begin{subfigure}{1.0\linewidth}
        \centering
         \includegraphics[width=1.0\columnwidth]{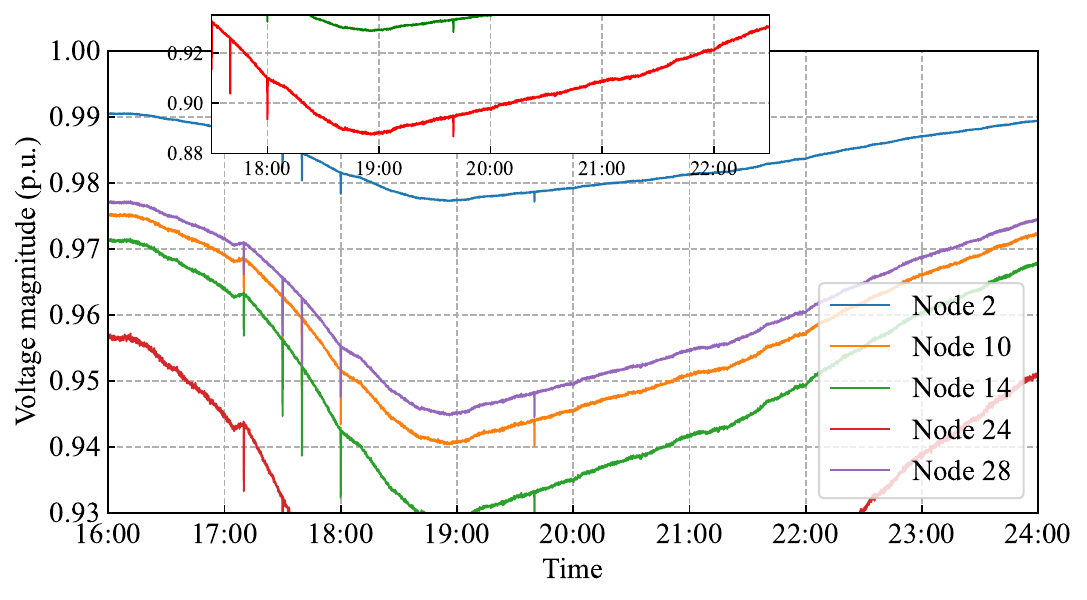}
         \caption{}
         \label{fig:voltages:nocontrol}
    \end{subfigure}
    \centering
    \begin{subfigure}{1.0\linewidth}
        \centering
         \includegraphics[width=1.0\columnwidth]{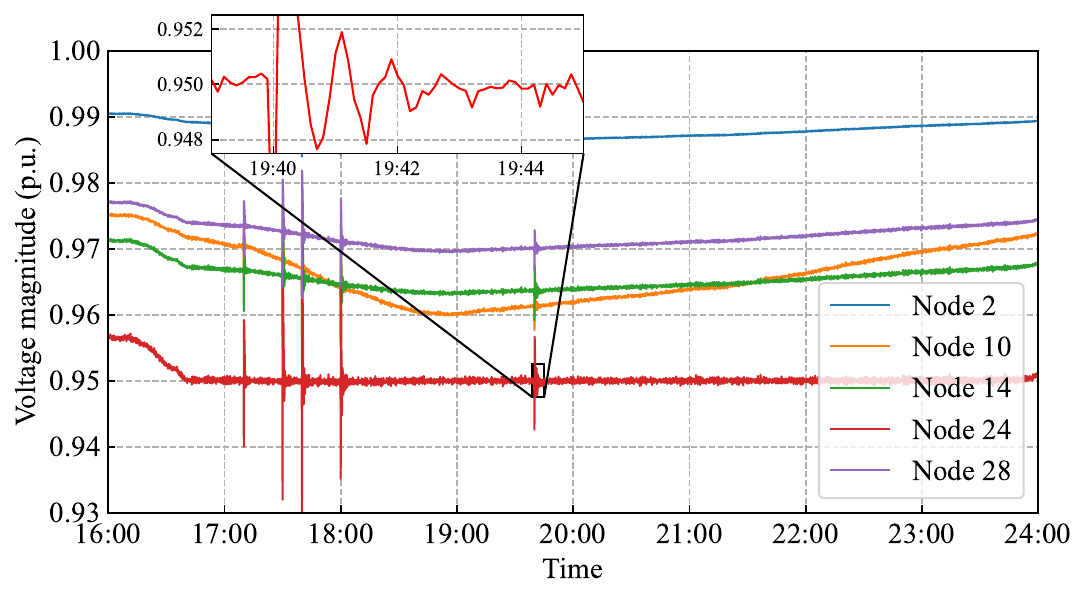}
         \caption{}
         \label{fig:voltages:PD}
    \end{subfigure}
    \centering
    \begin{subfigure}{1.0\linewidth}
        \centering
         \includegraphics[width=1.0\columnwidth]{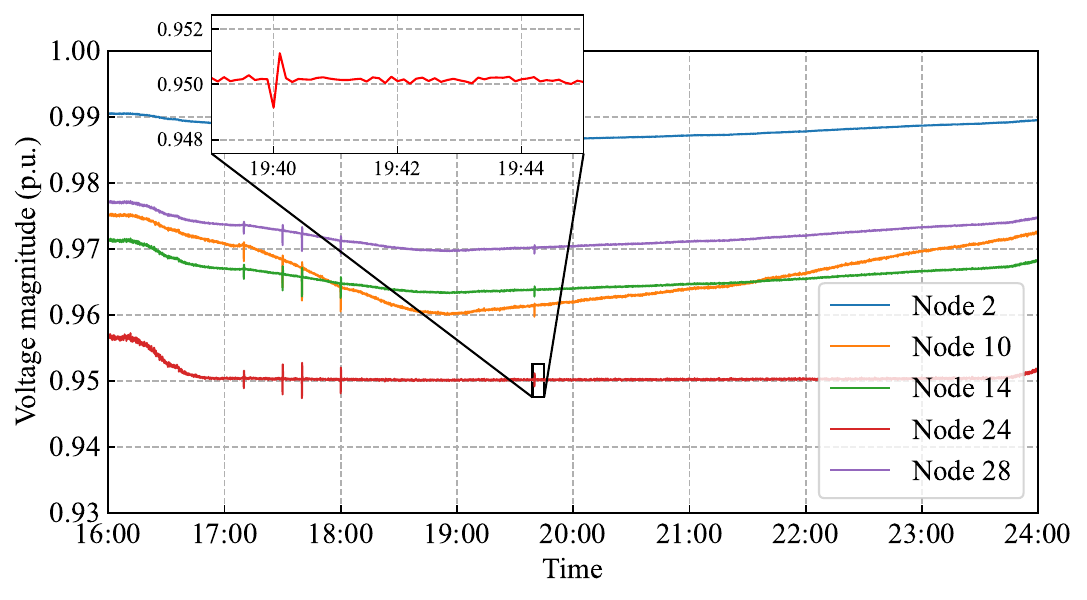}
         \caption{}
         \label{fig:voltages:control}
    \end{subfigure}
    \caption{The voltage profiles: (a) without control; (b) the primal-dual controller \cite{anese2018optimal}; (c) the proposed data-driven controller (trained with $\beta=0.1$).}
    \label{fig:voltages}
\end{figure}

We train the proposed data-driven control with the gradient calculated in \eqref{grad::ellphi}, and compare it with the primal-dual control from \cite{anese2018optimal} in the test set. We emphasize that the primal-dual control needs to transmit/receive signals to/from a central operator through a communication network, while the proposed control uses only the local measurements. As the solar generation declines and the net demand surges, if no action is taken, the feeder would undergo severe under-voltages as shown in Figure \ref{fig:voltages}(a). To enhance readability, only the voltage profiles at nodes \{2,10,14,24,28\} are provided. Some spikes are observed due to sudden variations of the solar and wind generations. Figure \ref{fig:voltages}(b) illustrates the voltage profiles under the primal-dual control. It regulates the voltages around the lower limit $0.95$ p.u. for most of the time, while experiencing heavy oscillations following the generation spikes. Figure \ref{fig:voltages}(c) shows the voltage profiles achieved by the proposed data-driven controller trained with the chance constraint threshold $\beta=0.1$. Not only does it enforce safe voltages, but it also quickly suppresses the flickers and drives the voltages above the lower limit in a few seconds. This is mainly because the proposed control updates the next operating point based on the feedback right before the current update, while the primal-dual control utilizes the dual variables accumulated over time and thus becomes less responsive. 

\begin{figure}
\centering
    \begin{subfigure}{1.0\linewidth}
        \centering
         \includegraphics[width=1.0\columnwidth]{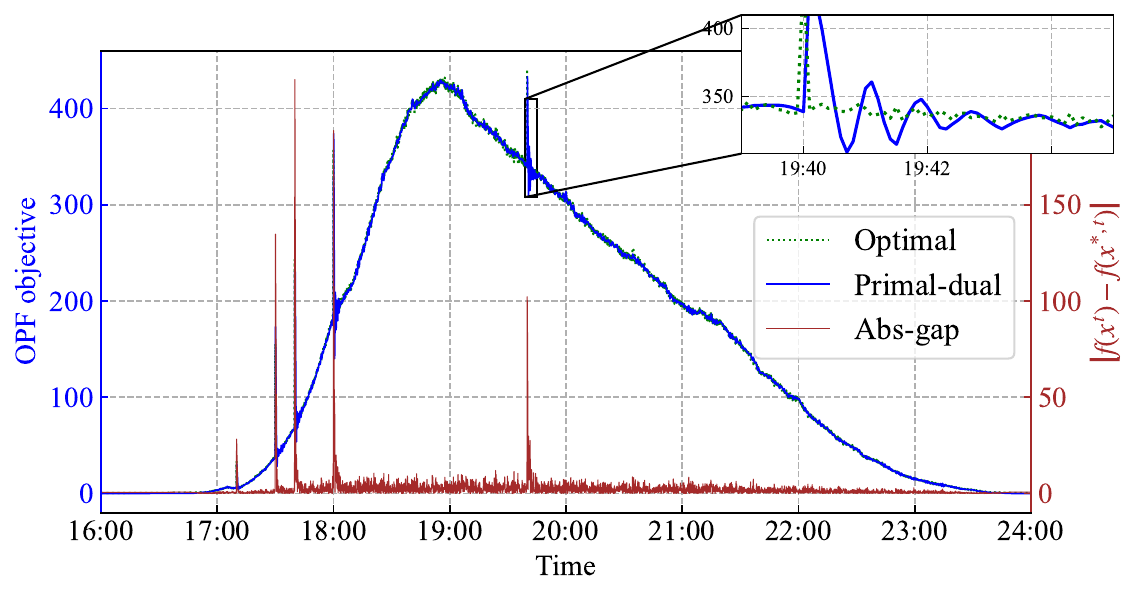}
         \caption{}
         \label{fig:tracking:PD}
    \end{subfigure}
    \centering
    \begin{subfigure}{1.0\linewidth}
        \centering
         \includegraphics[width=1.0\columnwidth]{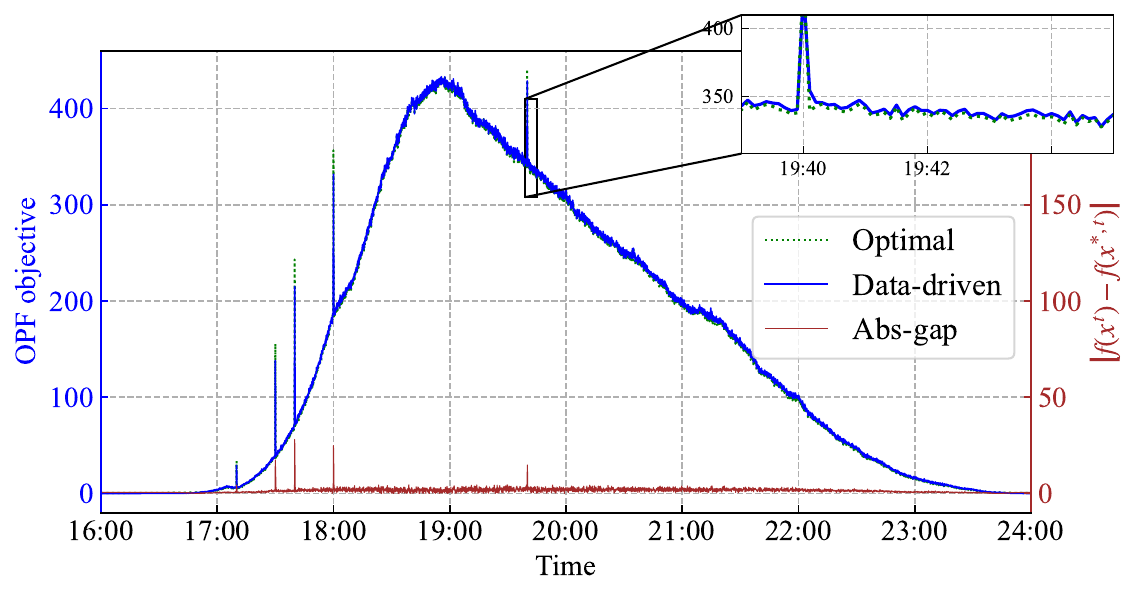}
         \caption{}
         \label{fig:tracking:control}
    \end{subfigure}
\caption{OPF objectives and the absolute gaps between the real-time operations $\boldsymbol{x}^t$ and optimal solutions $\boldsymbol{x}^{*,t}$: (a) the primal-dual controller \cite{anese2018optimal}; (b) the proposed data-driven controller (trained with $\beta=0.1$).}
\label{fig:obj_and_tracking}
\end{figure}

To further highlight the advantages of the proposed method, Figure \ref{fig:obj_and_tracking} compares the OPF objectives of the primal-dual control and the proposed control. The solid blue lines show the objective values obtained by different controls, and the dashed green lines are the (local) optimal values. The brown lines at the bottom illustrate the absolute gaps $\left|\boldsymbol{f}^t(\boldsymbol{x}^t)-\boldsymbol{f}^{t}(\boldsymbol{x}^{*,t})\right|$ between the controlled objective values and the (local) optimal values. It verifies that the proposed method can reduce the suboptimality gaps for most of the time. Especially near the generation spikes, the proposed method can quickly suppress the spikes and damp the oscillations.

We use the statistics below to quantify the performance:
\begin{subequations}
\begin{alignat}{2}
\text{Absolute-gap}&=\frac{1}{T}\sum_{t=1}^{T}\left|\boldsymbol{f}^t(\boldsymbol{x}^t)-\boldsymbol{f}^{t}(\boldsymbol{x}^{*,t})\right|, \nonumber\\
\text{Relative-gap}&=\frac{1}{T}\sum_{t=1}^{T}\frac{\left|\boldsymbol{f}^t(\boldsymbol{x}^t)-\boldsymbol{f}^{t}(\boldsymbol{x}^{*,t})\right|}{\boldsymbol{f}^{t}(\boldsymbol{x}^{*,t})},\nonumber\\
\text{Volt-violation}&=\frac{1}{T}\sum_{t=1}^{T}\left\|[\underline{\boldsymbol{V}}-\boldsymbol{V}^t]_{+}\right\|_2+\left\|[\boldsymbol{V}^t-\overline{\boldsymbol{V}}]_+\right\|_2. \nonumber
\end{alignat}
\end{subequations}
Table \ref{tab:comparePD} reports the above statistics of the primal-dual control and the proposed data-driven control. The proposed control is trained under different chance constraint thresholds $\beta=\{0.05,0.1,0.5\}$. We can see that with all $\beta$ values, the proposed control has a smaller absolute gap and relative gap, as well as a much smaller voltage violation, than the primal-dual control. Interestingly, the voltage violation decreases, while both the absolute gap and relative gap increase, as $\beta$ decreases. This is because a smaller $\beta$ reduces the tolerance on voltage violation and thus shrinks the feasible set, which makes the OPF solutions more conservative by pushing them further into the safe voltage region. This result also suggests a tradeoff between voltage safety and cost effectiveness, which can be controlled by tuning $\beta$ according to the practical need of the operator. 

\begin{table}
\centering
\caption{Comparison of the proposed data-driven control trained under different $\beta$, with the primal-dual control.}
\renewcommand\arraystretch{1.2}
\setlength\tabcolsep{4.5pt}
\begin{tabular}{c|c|c|c|c}
\toprule\hline
\multirow{2}{*}{\text {Metrics}} & \multirow{2}{*}{\makecell[c]{Primal-dual \\ control \cite{anese2018optimal}}} & \multicolumn{3}{c}{\text {Data-driven control}}   \\
\cline { 3 - 5 } &  & $\beta=0.05$ & $\beta=0.1$ & $\beta=0.5$ \\
\hline 
Absolute-gap & 2.248  & 1.961 & 1.466  & 0.948  \\
Relative-gap & 0.0173  & 0.0154 & 0.0113  & 0.0066   \\
Volt-violation & $3.5\times 10^{-4}$ & $6.4\times 10^{-6}$ & $6.9\times 10^{-6}$ & $3.5\times 10^{-5}$    \\
\hline\bottomrule
\end{tabular}
\label{tab:comparePD}
\end{table}
Besides, the proposed data-driven control is computationally efficient. The average computation time of each real-time operating point updated by the proposed control is $0.0059$ seconds, while that by the primal-dual control is $0.0083$ seconds. The proposed control is faster by $28.9\%$ than the primal-dual control, and the acceleration by the proposed control can be more significant in large-scale networks, owing to the use of local feedback that avoids the computation governing the entire network with its overhead growing quadratically with the network size.

\begin{figure}[t]
\centering
    \begin{subfigure}{0.49\linewidth}
        \centering
         \includegraphics[width=1.0\columnwidth]{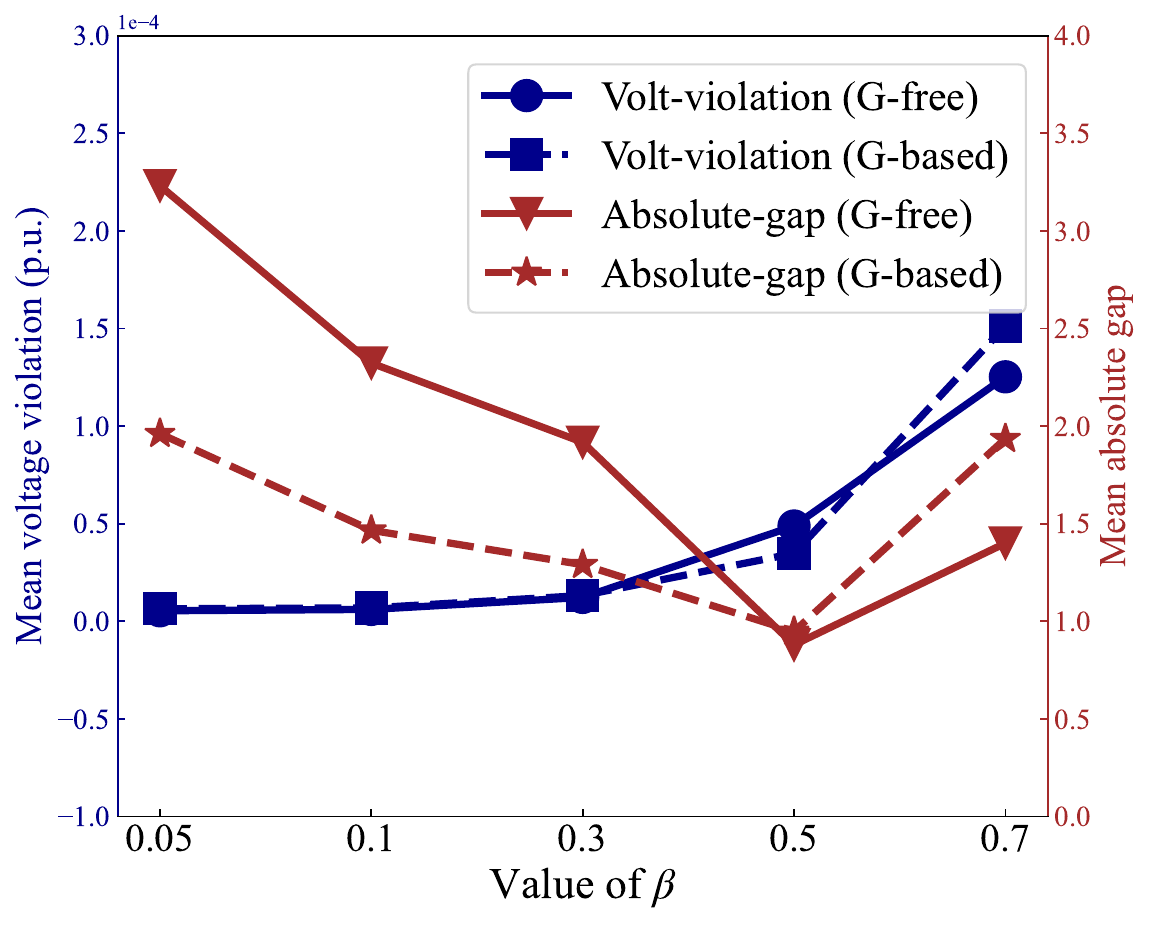}
         \caption{}
         \label{fig:trend_alpha}
    \end{subfigure}
    \centering
    \begin{subfigure}{0.49\linewidth}
        \centering
         \includegraphics[width=1.0\columnwidth]{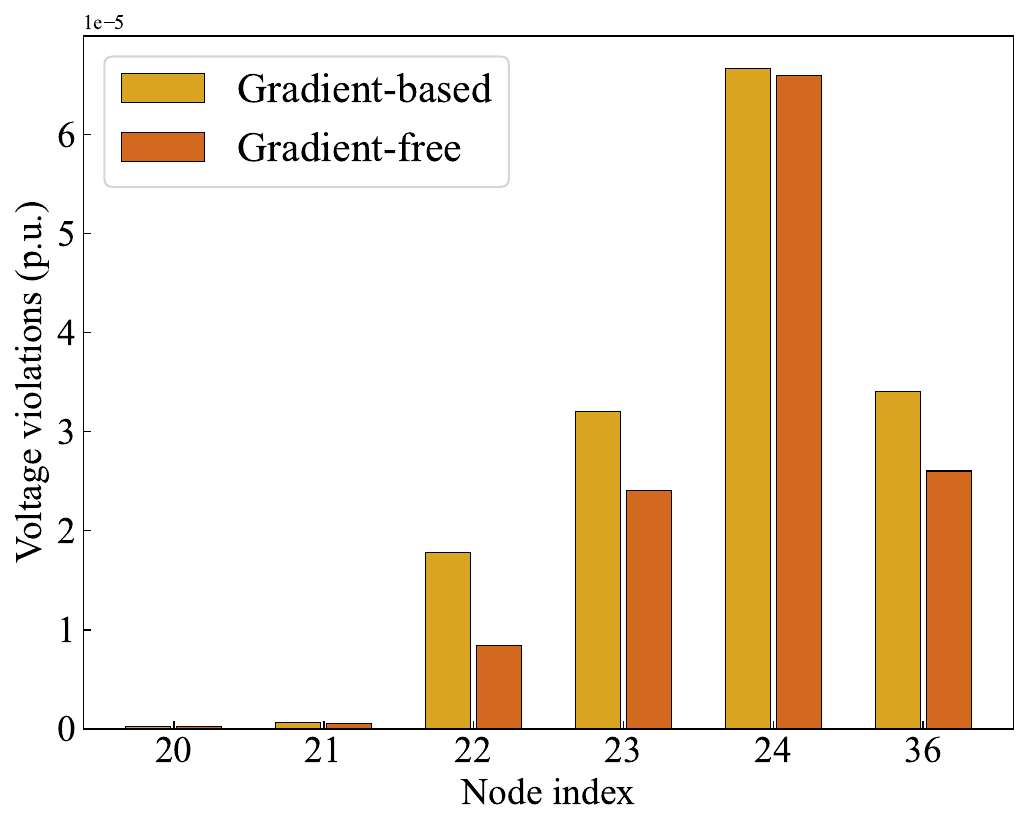}
         \caption{}
         \label{fig:basedvsfree}
    \end{subfigure}
\caption{Comparison of gradient-free learning and gradient-based learning: (a) the changes of average absolute gap and average voltage violation with different $\beta$; (b) the average voltage violations of nodes whose voltages have ever violated the safety limits under the training with $\beta=0.7$.}
\label{fig:gradientfree}
\end{figure}

\subsection{Gradient-Free Learning}
In this part, we train the proposed data-driven control in a gradient-free manner, in which the $2n$-point gradient estimator \eqref{grad::free} is applied to the nonlinear power flow model. The scale of the perturbation is set as $\varepsilon=0.001$. Figure \ref{fig:gradientfree}(a) shows the changes of the absolute gap in objectives and the voltage violation of the proposed control trained with five different values of $\beta=\{0.05,0.1,0.3,0.5,0.7\}$. The gradient-free learning approach is compared with its gradient-based counterpart. The solid lines show the results of the gradient-free training, while the dashed lines represent those of the gradient-based training. In general, the gradient-free approach can achieve comparable OPF objective values and voltage safety to the gradient-based approach, although the former avoids the complicated process of explicitly taking the gradient of the nonlinear power flow model. For the particular case $\beta=0.7$ in which the chance constraints for voltage safety are relatively relaxed, the gradient-free approach tends to reduce the voltage violation compared to the gradient-based approach, as further verified by the comparison of voltage violations at selected nodes in Figure \ref{fig:gradientfree}(b).

\begin{figure}[t]
    \centering
    \includegraphics[width=1.0\columnwidth]{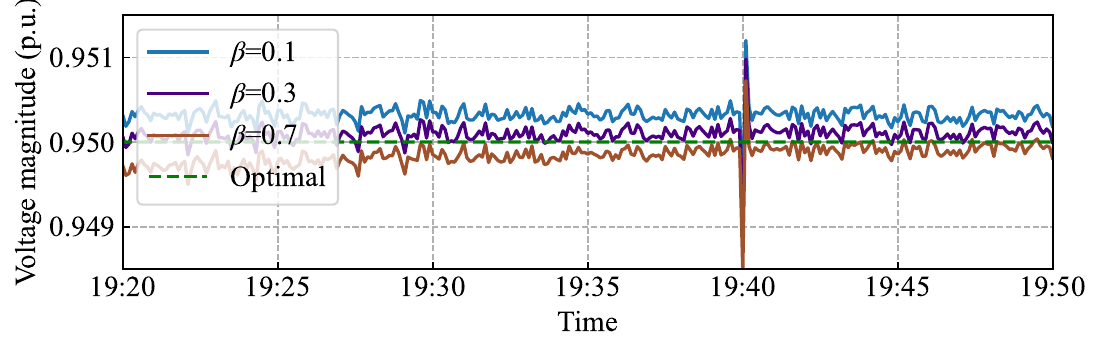}
    \caption{The voltage profiles of node $24$ when the proposed control is trained in a gradient-free manner with different $\beta$.}
    \label{fig:compare_alpha}
\end{figure}

It is also noted that for both gradient-free and gradient-based approaches, the absolute gap between the controlled objective values and the optimal values first decreases with $\beta$ until $\beta=0.5$ (as shown and explained in Section \ref{subsec::based}), and then increases with $\beta$. Indeed, as $\beta$ increases to relax the chance constraints for voltage safety, the proposed OPF solver tends to reduce the generations and thus the generation cost, which at the same time lowers the voltages. For very large $\beta$ such as $0.7$, the generation cost is even lower than the optimal value of OPF problem \eqref{OPF}, which is why the \emph{absolute} gap between them increases. However, this lower-cost solution is not feasible because it already violates voltage safety. This is verified in Figure \ref{fig:compare_alpha}, where the voltages obtained by the proposed control are lowered as $\beta$ increases. Particularly for $\beta=0.7$, the voltages drop even below those at the optimal feasible solution of OPF problem \eqref{OPF}, with the latter exactly at the lower voltage limit.

\section{Conclusion}
\label{sec::conclude}
We proposed a data-driven real-time algorithm, endowed with modified DNNs, to track the time-varying OPF solutions. The proposed algorithm utilized the local measurements to provide feedback control-based OPF solution for power network operations, while bypassing the need for real-time communications. A chance constraint reformulation was put forward to account for the changes and uncertainties of the time-varying OPF problem, which enabled us to further develop a stochastic primal-dual learning strategy, as well as its gradient-free counterpart based on the nonlinear power flow, to train the proposed algorithm without resort to any labeled data. The OPF solution-tracking performance of the proposed algorithm was analyzed using the universal approximation capability of DNNs and verified by the numerical results. In the future, we shall incorporate the second-order information of the OPF problems into the data-driven methods to provide agile updates towards the fast-changing solutions. We are also interested in optimizing the transition cost (e.g., the energy loss during charging/discharging the energy storage devices) in time-varying OPF problems using data-driven methods.




%

\appendices

\section{Proof of Theorem \ref{thm::1}}
\label{append::A}
First, we show the dynamics \eqref{dynamic} has an equilibrium point.
Defining $g_t(\boldsymbol{x}^{t-1}):=\boldsymbol{x}^{t-1}-\alpha\left(\nabla \boldsymbol{f}^t(\boldsymbol{x}^{t-1})+\boldsymbol{u}_{\boldsymbol{\varphi}}(\hat{\boldsymbol{v}}^{t},\boldsymbol{d}^t)\right)$, we aim to prove $\boldsymbol{x}^t=\left[g_t(\boldsymbol{x}^{t-1})\right]_{\boldsymbol{\mathcal{Y}}^t}$ is a continuous function. Taking the derivative of $g_t(\boldsymbol{x}^{t-1})$ with respect to $\boldsymbol{x}^{t-1}$, we have:
\begin{alignat}{2}
\frac{\partial g_t(\boldsymbol{x}^{t-1})}{\partial \boldsymbol{x}^{t-1}}=\boldsymbol{I}-\alpha\left(\nabla^2_{\boldsymbol{x}\boldsymbol{x}}\boldsymbol{f}^t+\frac{\partial \boldsymbol{u}_{\boldsymbol{\vartheta}}}{\partial \hat{\boldsymbol{v}}^{t}}\frac{\partial \hat{\boldsymbol{v}}^{t}}{\partial \boldsymbol{x}^{t-1}}\right),
\end{alignat}
where we use the condition \textit{C1)} that $u_{\varphi_i}(v_i,d_i)=u_{\vartheta_i}(v_i)+u_{\phi_i}(d_i)$, and $\boldsymbol{I}\in\mathbb{R}^{2N\times 2N}$ is an identity matrix. Noting that $\hat{\boldsymbol{v}}^t=\mathcal{A}\boldsymbol{x}^{t-1}+\boldsymbol{v}_{env}^t$, we have:
\begin{alignat}{2}
   \left\|\frac{\partial \hat{\boldsymbol{v}}^{t}}{\partial \boldsymbol{x}^{t-1}}\right\|_2=\left\|
   \mathcal{A}\right\|_2<\infty.
\end{alignat}
Notice $L_{\boldsymbol{\vartheta}}=\max_{i\in\mathcal{N}} L_{\vartheta_i}$ from condition \textit{C3)}. Moreover, By inequality \eqref{f:smooth} in Assumption \ref{ass::1}, we can conclude $\left\|\nabla^2_{\boldsymbol{x}\boldsymbol{x}}\boldsymbol{f}^t(\boldsymbol{x})\right\|_2\leqslant \xi,~\forall\boldsymbol{x}\in\boldsymbol{\mathcal{Y}}^t,t\in\mathcal{T}$. We then have:
\begin{alignat}{2}\label{inftybound:1}
\left\|\frac{\partial g_t(\boldsymbol{x}^{t-1})}{\partial \boldsymbol{x}^{t-1}}\right\|_2\leqslant 1+\alpha (\xi+L_{\boldsymbol{\vartheta}}\left\|\mathcal{A}\right\|_2)<\infty.
\end{alignat}
On the other hand, for the convex compact set $\boldsymbol{\mathcal{Y}}^t$ in \eqref{region:1}, we can show that:
\begin{alignat}{2} \label{inftybound:2}
\left\|\frac{\partial \boldsymbol{x}^t}{\partial g_t(\boldsymbol{x}^{t-1})}\right\|_2=\left\|\frac{\partial \left[g_t(\boldsymbol{x}^{t-1})\right]_{\boldsymbol{\mathcal{Y}}^t}}{\partial g_t(\boldsymbol{x}^{t-1})}\right\|_2\leqslant 1.
\end{alignat}
By \eqref{inftybound:1} and \eqref{inftybound:2}, it is easy to conclude that $\boldsymbol{x}^t=\left[g_t(\boldsymbol{x}^{t-1})\right]_{\boldsymbol{\mathcal{Y}}^t}$ is continuous in $\boldsymbol{\mathcal{Y}}^t$. Since $\boldsymbol{\mathcal{Y}}^t$ is convex and compact, by using Brouwer's fixed point theorem \cite{border1985fixed}, the dynamics \eqref{dynamic} has a fixed point that satisfies \eqref{equilibrium}.

We then show that the equilibrium of dynamics \eqref{dynamic} is unique by contradiction. Suppose both $(\boldsymbol{x}^{\dagger},\boldsymbol{v}^{\dagger})$ and $(\boldsymbol{x}^{\ddagger},\boldsymbol{v}^{\ddagger})$ are equilibrium points of dynamics \eqref{dynamic} and $\boldsymbol{x}^{\dagger}\neq \boldsymbol{x}^{\ddagger}$. We substitute them into \eqref{equilibrium::x} and take the difference:
\begin{alignat}{2}\label{appex:difference}
   \left\|\boldsymbol{x}^{\dagger}-\boldsymbol{x}^{\ddagger}\right\|_2\leqslant \left\|\boldsymbol{x}^{\dagger}-\boldsymbol{x}^{\ddagger}-\alpha\left(\nabla \boldsymbol{f}^t(\boldsymbol{x}^{\dagger})-\nabla \boldsymbol{f}^t(\boldsymbol{x}^{\ddagger})\right.\right.\nonumber\\
\left.\left.   +\boldsymbol{u}_{\boldsymbol{\varphi}}(\boldsymbol{v}^{\dagger},\boldsymbol{d}^t)-\boldsymbol{u}_{\boldsymbol{\varphi}}(\boldsymbol{v}^{\ddagger},\boldsymbol{d}^t)\right)\right\|_2.
\end{alignat}
Define the diagonal matrix $\boldsymbol{D}_{\boldsymbol{p}}\in\mathbb{R}^{N\times N}$ with its entries:
\begin{alignat}{2}\label{def::D}
\boldsymbol{D}_{p,i}=\left\{\begin{array}{cc}
\frac{u_{\vartheta_{p,i}}(v_i^{\dagger})-u_{\vartheta_{p,i}}(v_i^{\ddagger})}{v_i^{\dagger}-v_i^{\ddagger}} & v_i^{\dagger} \neq v_i^{\ddagger}, \\
0 & v_i^{\dagger}=v_i^{\ddagger}.
\end{array}\right.
\end{alignat}
The diagonal matrix $\boldsymbol{D}_{\boldsymbol{q}}\in\mathbb{R}^{N\times N}$ is defined similarly, and let $\boldsymbol{D}:=[\boldsymbol{D}_{\boldsymbol{p}}^{\top},\boldsymbol{D}_{\boldsymbol{q}}^{\top}]^{\top}$. By conditions \textit{C1)} and \textit{C2)}, we have:
\begin{alignat}{2}\label{appex:varphiv}
\boldsymbol{u}_{\boldsymbol{\varphi}}(\boldsymbol{v}^{\dagger},\boldsymbol{d}^t)-\boldsymbol{u}_{\boldsymbol{\varphi}}(\boldsymbol{v}^{\ddagger},\boldsymbol{d}^t)=\boldsymbol{D}\mathcal{A}(\boldsymbol{x}^{\dagger}-\boldsymbol{x}^{\ddagger}).
\end{alignat}
Substituting \eqref{appex:varphiv} into the RHS of \eqref{appex:difference}, we have:
\begin{alignat}{2}\label{appex:inequal}
&\left\|\boldsymbol{x}^{\dagger}-\boldsymbol{x}^{\ddagger}\right\|_2\leqslant\left\|\left(\boldsymbol{x}^{\dagger}-\boldsymbol{x}^{\ddagger}\right)-\alpha\left(\nabla\boldsymbol{f}^t(\boldsymbol{x}^{\dagger})-\nabla\boldsymbol{f}^t(\boldsymbol{x}^{\ddagger})\right)\right\|_2 \nonumber \\
&\qquad+\alpha\left\|\boldsymbol{D}\right\|_2\left\|\mathcal{A}\right\|_2 \left\|\boldsymbol{x}^{\dagger}-\boldsymbol{x}^{\ddagger}\right\|_2 \nonumber\\
&\leqslant (\sqrt{1-2\alpha m+\alpha^2\xi^2}+\alpha\left\|\boldsymbol{D}\right\|_2\left\|\mathcal{A}\right\|_2)\left\|\boldsymbol{x}^{\dagger}-\boldsymbol{x}^{\ddagger}\right\|_2,
\end{alignat}
where we use the properties \eqref{f:convex} and \eqref{f:smooth} from Assumption \ref{ass::1} in the second inequality. Recall condition \textit{C3)} that $L_{\boldsymbol{\vartheta}}=\max_{i\in\mathcal{N}}L_{\vartheta_i}<\frac{1-\sqrt{1-2\alpha m+\alpha^2\xi^2}}{\alpha \left\|\mathcal{A}\right\|_2}$, and notice the definition of matrix $\boldsymbol{D}$, which implies $\left\|\boldsymbol{D}\right\|_2< \frac{1-\sqrt{1-2\alpha m+\alpha^2\xi^2}}{\alpha \left\|\mathcal{A}\right\|_2}$. We have:
$$\sqrt{1-2\alpha m+\alpha^2\xi^2}+\alpha\left\|\boldsymbol{D}\right\|_2\left\|\mathcal{A}\right\|_2<1.$$
This implies the inequality \eqref{appex:inequal} holds only when $\boldsymbol{x}^{\dagger}-\boldsymbol{x}^{\ddagger}=0$, which leads to a contradiction. This completes the proof.

\section{Proof of Theorem \ref{thm::2}}
\label{append::B}

Consider the $\boldsymbol{x}^t$ generated by dynamics \eqref{dynamic}, equivalently \eqref{dynamic:compact}, at time $t$. We notice that:
\begin{alignat}{2}\label{inequality::xt}
    \left\|\boldsymbol{x}^t-\boldsymbol{x}^{*,t}\right\|_2\leqslant \left\|\boldsymbol{x}^t-\boldsymbol{x}^{\dagger,t}\right\|_2+\left\|\boldsymbol{x}^{\dagger,t}-\boldsymbol{x}^{*,t}\right\|_2,
\end{alignat}
where $\boldsymbol{x}^{\dagger,t}$ denotes the equilibrium point of \eqref{dynamic:compact} at time $t$. We first show the bound of the first term. Take the difference between dynamics \eqref{dynamic:compact} and the definition \eqref{equilibrium} of the equilibrium point. To differentiate $\hat{\boldsymbol{v}}^t$ in \eqref{dynamic:compact} and $\boldsymbol{v}^{\dagger,t}$ in \eqref{equilibrium}, we rewrite $\hat{\boldsymbol{v}}^t=\boldsymbol{v}^{t}(\boldsymbol{x}^{t-1})$ and $\boldsymbol{v}^{\dagger,t}=\boldsymbol{v}^{t}(\boldsymbol{x}^{\dagger,t})$ and have:
\begin{alignat}{2}\label{inequality::1}
    &\left\|\boldsymbol{x}^t-\boldsymbol{x}^{\dagger,t}\right\|_2^2\leqslant\left\|\boldsymbol{x}^{t-1}-\boldsymbol{x}^{\dagger,t}-\alpha[\nabla \boldsymbol{f}^t(\boldsymbol{x}^{t-1})-\nabla \boldsymbol{f}^t(\boldsymbol{x}^{\dagger,t})]\right.\nonumber\\
    &\left.\quad-\alpha\left[\boldsymbol{u}_{\boldsymbol{\varphi}}\left(\boldsymbol{v}^t(\boldsymbol{x}^{t-1}),\boldsymbol{d}^t\right)-\boldsymbol{u}_{\boldsymbol{\varphi}}\left(\boldsymbol{v}^t(\boldsymbol{x}^{\dagger,t}),\boldsymbol{d}^t\right)\right]\right\|_2^2 \nonumber\\
    &\leqslant\left\|\boldsymbol{x}^{t-1}-\boldsymbol{x}^{\dagger,t}\right\|_2^2+\alpha^2\left\|\nabla\boldsymbol{f}^t(\boldsymbol{x}^{t-1})-\nabla\boldsymbol{f}^t(\boldsymbol{x}^{\dagger,t})\right\|_2^2\nonumber\\
    &\quad+\alpha^2\left\|\boldsymbol{u}_{\vartheta}\left(\boldsymbol{v}^{t}({\boldsymbol{x}}^{t-1})\right)-\boldsymbol{u}_{\vartheta}\left(\boldsymbol{v}^{t}(\boldsymbol{x}^{\dagger,t})\right)\right\|_2^2 \nonumber\\
    &\quad-2\alpha\left(\boldsymbol{x}^{t-1}-\boldsymbol{x}^{\dagger,t}\right)^{\top}\left(\nabla\boldsymbol{f}^t(\boldsymbol{x}^{t-1})-\nabla\boldsymbol{f}^t(\boldsymbol{x}^{\dagger,t})\right)\nonumber\\
    &\quad -2\alpha\left(\boldsymbol{x}^{t-1}-\boldsymbol{x}^{\dagger,t}\right)^{\top}\left(\boldsymbol{u}_{\vartheta}\left(\boldsymbol{v}^{t}({\boldsymbol{x}}^{t-1})\right)-\boldsymbol{u}_{\vartheta}\left(\boldsymbol{v}^{t}(\boldsymbol{x}^{\dagger,t})\right)\right)\nonumber\\
     &\quad+2\alpha^2\left(\nabla\boldsymbol{f}^t(\boldsymbol{x}^{t-1})-\nabla\boldsymbol{f}^t(\boldsymbol{x}^{\dagger,t})\right)^{\top}\cdot\nonumber\\
    &\quad\quad\left(\boldsymbol{u}_{\vartheta}\left(\boldsymbol{v}^{t}({\boldsymbol{x}}^{t-1})\right)-\boldsymbol{u}_{\vartheta}\left(\boldsymbol{v}^{t}(\boldsymbol{x}^{\dagger,t})\right)\right),
\end{alignat}
where the first inequality uses the non-expansiveness of projection. The second inequality utilizes the property that $\boldsymbol{u}_{\boldsymbol{\varphi}}$ is separable over the inputs $\boldsymbol{v}^t$ and $\boldsymbol{d}^t$, i.e., $\boldsymbol{u}_{\boldsymbol{\varphi}}(\boldsymbol{v}^t,\boldsymbol{d}^t)=\boldsymbol{u}_{\boldsymbol{\vartheta}}(\boldsymbol{v}^t)+\boldsymbol{u}_{\boldsymbol{\phi}}(\boldsymbol{d}^t)$, thus eliminating $\boldsymbol{u}_{\boldsymbol{\phi}}(\boldsymbol{d}^t)$ with the same input $\boldsymbol{d}^t$. Noticing that $\boldsymbol{u}_{\boldsymbol{\vartheta}}$ is non-decreasing in $\boldsymbol{v}^t$ and Lipschitz continuous with constant $L_{\boldsymbol{\vartheta}}$, we have:
\begin{alignat}{2}\label{smooth::u}
&\left\|\boldsymbol{u}_{\vartheta}\left(\boldsymbol{v}^{t}({\boldsymbol{x}}^{t-1})\right)-\boldsymbol{u}_{\vartheta}\left(\boldsymbol{v}^{t}(\boldsymbol{x}^{\dagger,t})\right)\right\|_2 \nonumber\\ 
&\qquad\qquad\qquad \qquad  \leqslant L_{\boldsymbol{\vartheta}}\left\|\mathcal{A}\right\|_2\left\|\boldsymbol{x}^{t-1}-\boldsymbol{x}^{\dagger,t}\right\|_2,\nonumber\\
&(\boldsymbol{x}^{t-1}-\boldsymbol{x}^{\dagger,t})^{\top}\left(\boldsymbol{u}_{\vartheta}\left(\boldsymbol{v}^{t}({\boldsymbol{x}}^{t-1})\right)-\boldsymbol{u}_{\vartheta}\left(\boldsymbol{v}^{t}(\boldsymbol{x}^{\dagger,t})\right)\right)\geqslant 0.
\end{alignat}

Substituting inequalities \eqref{fproperty}, \eqref{smooth::u} into the right-hand-side (RHS) of \eqref{inequality::1}, we have:
\begin{alignat}{2}
    \left\|\boldsymbol{x}^t-\boldsymbol{x}^{\dagger,t}\right\|_2^2&\leqslant \big(1+\alpha^2\xi^2+\alpha^2 L_{\boldsymbol{\vartheta}}^2\left\|\mathcal{A}\right\|_2^2+2\alpha^2\xi L_{\boldsymbol{\vartheta}}\left\|\mathcal{A}\right\|_2\nonumber \\
    &\qquad-2\alpha m \big)\left\|\boldsymbol{x}^{t-1}-\boldsymbol{x}^{\dagger,t}\right\|^2_2\nonumber\\
    &=\rho^2(\alpha)\left\|\boldsymbol{x}^{t-1}-\boldsymbol{x}^{\dagger,t}\right\|_2^2, \nonumber
\end{alignat}
i.e.,
\begin{alignat}{2}\label{inequality::t_1}
    \left\|\boldsymbol{x}^t-\boldsymbol{x}^{\dagger,t}\right\|_2\leqslant \rho(\alpha)\left\|\boldsymbol{x}^{t-1}-\boldsymbol{x}^{\dagger,t}\right\|_2,
\end{alignat}
by the definition of $\rho(\alpha)$. Substituting \eqref{inequality::t_1} into \eqref{inequality::xt}, we have:
\begin{alignat}{2} \label{inequality::star}
    &\quad\left\|\boldsymbol{x}^t-\boldsymbol{x}^{*,t}\right\|_2\leqslant \rho(\alpha)\left\|\boldsymbol{x}^{t-1}-\boldsymbol{x}^{\dagger,t}\right\|_2+\left\|\boldsymbol{x}^{\dagger,t}-\boldsymbol{x}^{*,t}\right\|_2 \nonumber\\
    &\leqslant \rho(\alpha)\left(\left\|\boldsymbol{x}^{t-1}-\boldsymbol{x}^{*,t}\right\|_2+\left\|\boldsymbol{x}^{*,t}-\boldsymbol{x}^{\dagger,t}\right\|_2\right)+\left\|\boldsymbol{x}^{\dagger,t}-\boldsymbol{x}^{*,t}\right\|_2\nonumber\\
    &=\rho(\alpha)\left\|\boldsymbol{x}^{t-1}-\boldsymbol{x}^{*,t}\right\|_2+\left(1+\rho(\alpha)\right)\left\|\boldsymbol{x}^{\dagger,t}-\boldsymbol{x}^{*,t}\right\|_2.
\end{alignat}

We next derive the bound to the second term in \eqref{inequality::xt}. 
Before proceeding, we prove Lemma \ref{lemm::1}. Let $g_t(\boldsymbol{x}^{\dagger,t})=\boldsymbol{x}^{\dagger,t}-\alpha\left(\nabla\boldsymbol{f}^t(\boldsymbol{x}^{\dagger,t})+\boldsymbol{u}_{\boldsymbol{\varphi}}(\boldsymbol{v}^{\dagger,t},\boldsymbol{d})\right)$, and take the derivative on both sides of \eqref{Ht} with respect to $\boldsymbol{u}_{\boldsymbol{\varphi}}$, we have:
\begin{alignat}{2}
   \left\| \frac{\partial H_t(\boldsymbol{\varphi})}{\partial \boldsymbol{u}_{\boldsymbol{\varphi}}}\right\|_2=\left\|\frac{\partial \left[ g_t(\boldsymbol{x}^{\dagger,t})\right]_{\boldsymbol{\mathcal{Y}^t}}}{\partial g_t(\boldsymbol{x}^{\dagger,t})}\cdot\frac{\partial g_t(\boldsymbol{x}^{\dagger,t})}{\partial \boldsymbol{u}_{\boldsymbol{\varphi}}}\right\|_2\leqslant \alpha, \nonumber
\end{alignat}
where we use the result from \eqref{inftybound:2}, and $\left\|\partial g_t(\boldsymbol{x}^{\dagger,t})/\partial \boldsymbol{u}_{\boldsymbol{\varphi}}\right\|_2\leqslant \alpha$. This proves Lemma \ref{lemm::1}. We further consider:
\begin{alignat}{2}
\boldsymbol{x}^{\dagger,t}=H_t(\boldsymbol{\varphi^*}), \text{and } \boldsymbol{x}^{*,t}=H_t(\boldsymbol{\varphi}^{*,t}),
\end{alignat}
where $\boldsymbol{\varphi}^*$ denotes the parameter returned by \eqref{primaldual} that generates $\boldsymbol{x}^{\dagger,t}$, and $\boldsymbol{\varphi}^{*,t}$ generates $\boldsymbol{x}^{*,t}$  by solving \eqref{OPFphi} separately at each $t$. Using Lemma \ref{lemm::1}, we have:\begin{alignat}{2}\label{inequality::dagger}
    &\left\|\boldsymbol{x}^{\dagger,t}-\boldsymbol{x}^{*,t}\right\|_2=\left\|H_t(\boldsymbol{\varphi}^{*})-H_t(\boldsymbol{\varphi}^{*,t})\right\|_2\nonumber\\
    &\leqslant L_{h}\left\|\boldsymbol{u}_{\boldsymbol{\varphi}^{*}}-\boldsymbol{u}_{\boldsymbol{\varphi}^{*,t}}\right\|_2\leqslant L_{h}\epsilon,
\end{alignat}
where the last inequality comes from the $\epsilon$-universal property of $\boldsymbol{u}_{\boldsymbol{\varphi}^{*}}$. Besides, by triangle inequality and \eqref{xstarchange}, we also have:
\begin{alignat}{2} \label{x::t_1_star}
    \left\|\boldsymbol{x}^{t-1}-\boldsymbol{x}^{*,t}\right\|_2&=\left\|\boldsymbol{x}^{t-1}-\boldsymbol{x}^{*,t-1}+\boldsymbol{x}^{*,t-1}-\boldsymbol{x}^{*,t}\right\|_2\nonumber\\
    &\leqslant\left\|\boldsymbol{x}^{t-1}-\boldsymbol{x}^{*,t-1}\right\|_2+\gamma.
\end{alignat}

Substituting \eqref{inequality::dagger} and \eqref{x::t_1_star} in to the RHS of \eqref{inequality::star}, we have:
\begin{alignat}{2}
    &\left\|\boldsymbol{x}^t-\boldsymbol{x}^{*,t}\right\|_2\leqslant \rho(\alpha)\left\|\boldsymbol{x}^{t-1}-\boldsymbol{x}^{*,t-1}\right\|_2\nonumber\\
    &\qquad+\rho(\alpha)\gamma+\left(1+\rho(\alpha)\right)L_{h}\epsilon \nonumber\\
    & \leqslant\rho^t(\alpha) \left\|\boldsymbol{x}^{0}-\boldsymbol{x}^{*,0}\right\|_2\nonumber\\
    &\qquad+\frac{1-\rho^t(\alpha)}{1-\rho(\alpha)}\big[\rho(\alpha)\gamma+\left(1+\rho(\alpha)\right)L_{h}\epsilon\big].\nonumber
\end{alignat}

With $\rho(\alpha)<1$, we can conclude that $\left\|\boldsymbol{x}^t-\boldsymbol{x}^{*,t}\right\|_2$ converges exponentially to the asymptotic bound: 
\begin{alignat}{2}
        \limsup_{t\rightarrow \infty}\left\|\boldsymbol{x}^{t}-\boldsymbol{x}^{*,t}\right\|_2=\frac{\rho(\alpha)\gamma+(1+\rho(\alpha))L_{h}\epsilon}{1-\rho(\alpha)},
\end{alignat}
which completes the proof.


 




\vfill

\end{document}